\documentclass[12pt,a4paper,english,reqno%,dvipdfmx
]{amsart}
\usepackage[utf8]{inputenc}
\usepackage{mathrsfs}
\usepackage{amsfonts}
\usepackage{amsthm}
\usepackage{amssymb}
\usepackage{amscd}
\usepackage{amstext}
\usepackage[all]{xy}
\xyoption{2cell}
\UseTwocells
\usepackage{graphicx}
\usepackage{comment}
\usepackage{calrsfs}%% \mathcal changes to \mathscr
\usepackage{framed}%%%% encycle by a frame

\usepackage[colorinlistoftodos]{todonotes}
\usepackage[colorlinks=true, allcolors=blue]{hyperref}

\usepackage{tikz}
\usetikzlibrary{positioning,intersections,cd,matrix, arrows, backgrounds}
\usetikzlibrary{shapes,shapes.geometric,shapes.misc}%% needed for rounded rectangle
\usetikzlibrary{decorations}
\usetikzlibrary{decorations.pathreplacing}%% instead of snakes
\newcommand\Ar[3]{\ar[from={#1}, to={#2}, #3]}
\newcommand\Nname[1]{|[alias=#1]|}

\usepackage[OT2,T1]{fontenc}

\theoremstyle{plain}
\newtheorem{thm}{Theorem}[section]
\newtheorem{prp}[thm]{Proposition}
\newtheorem{lem}[thm]{Lemma}
\newtheorem{cor}[thm]{Corollary}

\newtheorem{prb}[thm]{Problem}

\newtheorem*{thm-nn}{Theorem}
\newtheorem*{prp-nn}{Proposition}
\newtheorem*{lem-nn}{Lemma}
\newtheorem*{cor-nn}{Corollary}
\newtheorem*{clm-nn}{Claim}
\newtheorem*{cnj-nn}{Conjecture}
\newtheorem*{prb-nn}{Problem}

\theoremstyle{definition}
\newtheorem{dfn}[thm]{Definition}
\newtheorem{exm}[thm]{Example}

\newtheorem*{dfn-nn}{Definition}

\newtheorem{rmk}[thm]{Remark}
\newtheorem{ntn}[thm]{Notation}
\newtheorem{cvn}[thm]{Convention}

\newcommand{\xyR}[1]{%
\xydef@\xymatrixrowsep@{#1}}

\newcommand{\xyC}[1]{%
\xydef@\xymatrixcolsep@{#1}}

\newcommand\al{\alpha}
\newcommand\be{\beta}

\newcommand\de{\delta}
\newcommand\ep{\varepsilon}
\newcommand\ze{\zeta}
\newcommand\et{\eta}
\newcommand\ro{\rho}
\newcommand\si{\sigma}

\newcommand\up{\upsilon}

\newcommand\Ga{\Gamma}

\newcommand\La{\Lambda}

\newcommand\soc{\operatorname{soc}}
\newcommand\Ker{\operatorname{Ker}}
\newcommand\Cok{\operatorname{Coker}}
\renewcommand\Im{\operatorname{Im}}%%%%%%%%%%%%%%%%%%%%%%

\newcommand\Hom{\operatorname{Hom}}

\newcommand\rad{\operatorname{rad}}

\newcommand\End{\operatorname{End}}
\newcommand\Ext{\operatorname{Ext}}

\newcommand\rep{\operatorname{rep}}
\renewcommand\mod{\operatorname{mod}}%%%%%%%%%%%%%%%%%%%%%

\newcommand\supp{\operatorname{supp}}

\newcommand{\udim}{\operatorname{\underline{dim}}\nolimits}

\newcommand\prj{\operatorname{prj}}

\newcommand\add{\operatorname{add}}%-- added on 2008-07-29 ---

\newcommand\calC{{\mathcal C}}

\newcommand\calH{{\mathcal H}}
\newcommand\calI{{\mathcal I}}

\newcommand\calK{{\mathcal K}}
\newcommand\calL{{\mathcal L}}

\newcommand\calP{{\mathcal P}}

\newcommand\calS{{\mathcal S}}
\newcommand\calT{{\mathcal T}}

\newcommand\calX{{\mathcal X}}

\newcommand\bbA{{\mathbb A}}
\newcommand\bbI{{\mathbb I}}

\newcommand\bbZ{{\mathbb Z}}

\newcommand\bbR{{\mathbb R}}

\newcommand\op{^{\mathrm{op}}} % changed text --> mathrm

\newcommand\inv{^{-1}}

\renewcommand\implies{\text{$\Rightarrow$}\ }%%%%%%%%%%%%%%%%%%
\newcommand\equivalent{\text{$\Leftrightarrow$}\ }

\newcommand\iso{\cong}
\newcommand\ds{\oplus}
\newcommand\ox{\otimes}

\newcommand\ovl{\overline}
\newcommand\Ds{\bigoplus}

\def\dsm#1,#2..#3{\bigoplus_{{#1}={#2}}^{#3}}
\def\sm#1,#2..#3{\sum_{{#1}={#2}}^{#3}}

\newcommand\id{1\kern-.25em{\text{{\rm l}}}} %1\!\!{\text{{\rm l}}}   is good for 12pt

\newcommand\isoto{\ \raise.8ex\hbox{$^{\sim}$}\kern-.7em\hbox{$\to$}\ }
\newcommand\Cdot{\raisebox{2pt}{$\centerdot$}}
\newcommand\down{_{\Cdot}}
\renewcommand\up{^{\Cdot}}%%%%%%%%%%%%%%%%%%%%%%%%%%%%%%%%%%%%%%%%%%%%

\newcommand\ya[1]{\xrightarrow{#1}}
\newcommand\blank{\operatorname{-}}

\newcommand\bg{%
\family{cmr}\size{20}{12pt}\selectfont}

\newcommand\bigzerou{%
\smash{\lower1.7ex\hbox{\bg 0}}}

\def\repr[#1;#2;#3;#4;#5]{
\left(
\begin{matrix}#1\\#2\end{matrix}
#3
\begin{matrix}#4\\#5\end{matrix}
\right)}
\newcommand\sbmat[1]{\left[\begin{smallmatrix} #1 \end{smallmatrix}\right]}

\newcommand\bfL{\mathbf{L}}
\newcommand\bfP{\mathbf{P}}

\renewcommand\k{\Bbbk}%%%%%%%%%%%%%%%%%%%%%%%%%%%%%%%%%%

\newcommand\dom{\operatorname{dom}}
\newcommand\cod{\operatorname{cod}}

\newcommand\To{\Rightarrow}
\newcommand{\bi}[3]{{}_{#2}{#1}_{#3}}

\newcommand\dm{{}}
\newcommand{\Sint}{\mathrm{S}_{\mathrm{int}}}
\newcommand{\Fint}{\mathrm{F}_{\mathrm{int}}}

\newcommand{\cov}{\operatorname{Cov}}
\newcommand{\Seg}{\operatorname{Seg}}

\newcommand{\Mon}{\operatorname{Mon}}
\newcommand{\Epi}{\operatorname{Epi}}
\newcommand{\inj}{\operatorname{inj}}

\newcommand\adj{\dashv}
\newcommand{\comma}[2]{#1\!\downarrow\!#2}

\newcommand{\bcalK}{\bar{\calK}}

\begin{document}
\title[Koszul coresolutions]{Relative Koszul coresolutions and relative Betti numbers}

\author{Hideto Asashiba}

%Key words:
%Koszul complexes; resolutions; Betti numbers; interval modules; approximations

\makeatletter
\@namedef{subjclassname@2020}{%
  \textup{2020} Mathematics Subject Classification}
\makeatother
\thanks{The author is partially supported by JSPS Grant-in-Aid for Scientific Research
 (C) 18K03207, JSPS Grant-in-Aid for Transformative Research Areas (A) (22A201),
and by Osaka Central Advanced Mathematical Institute
(MEXT Promotion of Distinctive Joint Research Center Program
JPMXP0723833165%
%JPMXP0619217849
).}

\address{Department of Mathematics, Faculty of Science, Shizuoka University,
836 Ohya, Suruga-ku, Shizuoka, 422-8529, Japan;}
\address{
Institute for Advanced Study, KUIAS, Kyoto University,
Yoshida Ushinomiya-cho, Sakyo-ku, Kyoto 606-8501, 
Japan; and}
\address{Osaka Central Advanced Mathematical Institute,
3-3-138 Sugimoto, Sumiyoshi-ku,
Osaka, 558-8585, Japan.}
\email{asashiba.hideto@shizuoka.ac.jp}

\maketitle

\begin{abstract}
Let $G$ be a finitely generated
right $A$-module for a finite-di\-men\-sion\-al algebra $A$ over a fieled $\k$,
and $\calI$ the additive closure of $G$. 
We will define an $\calI$-relative Koszul coresolution $\calK\up(V)$
of an indecomposable direct summand $V$ of $G$, and
show that for a finitely generated $A$-module $M$,
the $\calI$-relative $i$-th Betti number for $M$ at $V$
is given as the $\k$-dimension of the $i$-th homology of the
$\calI$-relative Koszul complex $\calK_V(M)\down:=\Hom_A(\calK\up(V),M)$
of $M$ at $V$ for all $i \ge 0$.
This is applied to investigate the minimal interval resolution/coresolution
of a persistence module $M$,
e.g., to check the interval decomposability of  $M$, and to compute
the interval approximation of $M$.
\end{abstract}

\section{Introduction}

In topological data analysis (TDA for short),
persistent homology plays an important role in examining the topological
 properties 
 of the data \cite{MR1949898},
in which context the data are usually given in the form of a point cloud
(a finite subset of a finite-dimensional Euclidean space).  Given one-parameter filtrations
 of simplicial complexes 
 arising from the data,
the persistent homology  construction 
yields representations of a Dynkin quiver $Q$ of type $A$,  and  thus 
 to  modules over the path category of Q,
that are sometimes called 1-{\em  parameter  persistence modules} \cite{MR2121296,BL,MR3868218}.
The product quiver of $d$ Dynkin quiver of type $A$ with full commutativity relations for some $d \ge 1$ is called
$d$D-grid.
By considering multi-parameter filtrations of simplicial complexes,
representations of $d$D-grid naturally arise in practical settings,
which are called $d$-{ parameter} persistence modules \cite{MR3868218,MR3824276}.
Since the linear category defined by this quiver with relations
can be regarded as the incidence category
%(Definition \ref{dfn:inc-cat}(1)) 
of a poset,
persistence modules are understood as
modules over the incidence category  $\k[\bfP]$ of a poset $\bfP$
 over a field $\k$ in general,
or equivalently, functors from  $\bfP$ (regarded as a category) to the category
$\mod \k$ of finite-dimensional vector spaces over $\k$.
In this paper we restrict ourselves to finite posets.
In this case,
the incidence category  $\k[\bfP]$ can be regarded as the incidence algebra $\k\bfP$
(we often identify these),
and persistence modules are nothing but modules over $\k\bfP$.

When we discuss in general,
we consider a finite-dimensional algebra $A$ over $\k$ instead of $\k\bfP$, and
we assume that all $A$-modules in consideration are finite-dimensional over $\k$
unless otherwise stated.
We mainly work over right $A$-modules
(when $A = \k\bfP$ is regarded as a category, they are contravariant functors $\k\bfP \to \mod\k$),
and the category of right $A$-modules is denoted by $\mod A$.
The full subcategory of $\mod A$ consisting of the projective modules is denoted by $\prj A$.

\subsection{Relative Betti numbers}
\label{ssec:rel-Betti}
For a persistence module $M \in \mod \k\bfP$, recall that if
\begin{equation}
\label{eq:standard-resol}
\cdots \to \Ds_{a \in \bfP} P_a^{\be_M^1(a)} \to \Ds_{a \in \bfP} P_a^{\be_M^0(a)} \to M \to 0 
\end{equation}
is a minimal projective resolution of $M$,
where $P_a$ is the projective indecomposable $\k\bfP$-module corresponding to each
$a \in \bfP$ (i.e., as a functor $P_a = \k\bfP(\blank, a)$),
then the number $\be_M^i(a)$ is uniquely determined by $M$, and
is called the $i$-th \emph{Betti number} for $M$ at $a$ (\cite{lesnick2022computing}).
A minimal projective resolution of $M$ gives important (homological) information of $M$
using projective modules and morphisms between them,
and the Betti numbers give us all the terms of the resolution.
Thus Betti numbers are useful to investigate the structure of $M$.
\par
Now let $A$ be an algebra, $\{V_I \mid I \in \bbI\}$ a finite set of mutually non-isomorphic
indecomposable modules in $\mod A$,
and $\calI$ the full subcategory of $\mod A$ consisting of all finite direct sums
of modules $V_I,\ (I \in \bbI)$.
Analogous to projective resolutions of a module $M$ over an algebra $A$
in homological algebra as above,
there is a notion of $\calI$-resolutions of $M$ in $\calI$-relative
homological algebra, which measures $M$ by using modules
in $\calI$ instead of projective modules
(see e.g., \cite{auslander1993relative, enochs2011relative, hochschild1956relative},
also \cite{BBH, BBH2023exact} for applications to TDA).
Analogous to a minimal projective resolution, 
a minimal $\calI$-resolution
\[
\cdots \to \Ds_{I \in \calI} V_I^{\be_M^1(I)}
\to \Ds_{I \in \calI} V_I^{\be_M^0({ I})} \to M \to 0 
\]
of $M$ is defined by $\calI$-relative homological algebra $($Definition \ref{dfn:min-resol-Betti}$)$.
Similar to the case of the standard homological algebra,
the numbers $\be_M^i(a)$ are uniquely determined by $M$, and
is called the $\calI$-relative $i$-th \emph{Betti number} for $M$ at $a$ (\cite{botnan2022bottleneck}).

\subsection{Koszul complex}
\label{ssec:Kcomplex}
In the paper \cite{CGRST}, Chach{\'o}lski et.\,al.\ gave a way (Theorem 3.8) to compute
Betti numbers for a persistence module $M$ in the case where
$\bfP$ is an upper semi-lattice, which we now explain.
To distinguish it from the general poset, let us denote this $\bfP$ as $\bfL$.
Since we are mainly deals with right modules (i.e., modules over the opposite poset,
or contravariant functors $\k\bfL \to \mod \k$),
we translate the original construction into the right module setting below.
Hence in this explanation, we assume that $\bfL$ is a lower semi-lattice, and
set $\bfL = \{a_1, \dots, a_n\}$ ($a_i$ are pairwise distinct).

\begin{dfn}
Let $M$ be a right $\k\bfL$-module, and $a \in \bfL$.
Then the \emph{Koszul complex} $\calK_a(M)\down$ of $M$ at $a$
is a complex
\[
\cdots \to \calK_a(M)_d \ya{\partial} \cdots \ya{\partial} \calK_a(M)_1 \ya{\partial}\calK_a(M)_0 \to 0
\]
of vector spaces defined as follows:
$$
\calK_a(M)_0:= M(a),\, \calK_a(M)_i:= \Ds_{S \in \cov(a)_i}M(\vee S),\text{ for }i > 0,
$$
where
$
\cov(a):= \{b \in \bfL \mid a < b,\text{ and } a \le c \le b \text{ implies }c \in \{a,b\}\}
$
and
$
\cov(a)_i:= \{S \subseteq \cov(a) \mid \# S = i, S \text{ has an upper bound}\}%%%%%%%%%%%%%%%%%
\footnote{
Note here that since $\bfL$ is upper semi-lattice and $S$ has an upper bound,
$S$ has its join $\vee S$.},
$
and the differential $\partial_i \colon \Ds_{S \in \cov(a)_i}M(\vee S) \to
\Ds_{T \in \cov(a)_{i-1}}M(\vee T)\ (i > 0)$ is given by the matrix
$$
[\chi(T \subseteq S)(-1)^{n(S,T)}M(p_{\vee S, \vee T})]_{T, S},
$$
(replace $\vee T$ by $a$ if $i = 1$)
where 
$\chi(T \subseteq S):= 1$ if $T \subseteq S$, and $:= 0$ otherwise; and
$n(S,T):= j$ if $S = \{a_{s_0}, \cdots, a_{s_{i-1}}\}$ with $s_0 < \cdots < s_{i-1}$ and
$S \setminus T = \{a_{s_j}\}$; and $p_{b,a}$
denotes the unique morphism from $a$ to $b$ in $\bfL$ as a category
for all $a,\, b \in \bfL$.
\end{dfn}

Then they have shown the following (\cite[Theorem 3.8]{CGRST})\footnote{%%%%%%%%%%%%%%%
Acturally the authors replaced the lower semi-lattice assumption on $\bfL$
by the assumption that all $S \in \cov(a)_i$ has the join $\vee S$.}.

\begin{thm}
\label{thm:standardKoszul}
Let $M$ be a right $\k\bfL$-module, and $a \in \bfL$.
Then for each $d \ge 0$,
it holds that
\[
\be_M^d(a) = \dim H_d(\calK_a(M)\down).
\]
\end{thm}

We remark that in the paper \cite{CGRST}, the authors consider an $\k\bfL$-$A$-bimodule $\calT$
and relate the $\calI$-relative resolution in $\mod A$ with the standard projective resolution
\eqref{eq:standard-resol} in $\mod \k\bfL$
by using the adjoint pair $\blank\ox_{\k\bfL}\calT \adj \Hom_A(\calT,\blank)$, which makes it possible
to apply Theorem \ref{thm:standardKoszul} for $\mod \k\bfL$ to $\mod A$.
Instead, we here think that
it would be good to have a complex also in the relative version in $\mod A$
having the same property as the Koszul complex above.
This suggests us the following problem.
\begin{prb}
Let $A$, $M$, $\{V_I \mid I  \in \bbI\}$, and $\calI$ be as in Sect.\ \ref{ssec:rel-Betti}.
Let
\[
\cdots \to \Ds_{I \in \calI} V_I^{\be_M^1(a)} \to \Ds_{I \in \calI} V_I^{\be_M^0(a)} \to M \to 0 
\]
be a minimal $\calI$-resolution of $M$.
In this setting, define a complex $\calK_I(M)\down$ such that
\begin{equation}
\label{eq:Koszul-cpx-cdn}
\be_M^d(I) = \dim H_d(\calK_I(M)\down).
\end{equation}
\end{prb}

Analogous to the above, we call the defined complex $\calK_I(M)\down$
a \emph{Koszul complex} of $M$ at $I$.

\subsection{Sketch of our solution}
Keep the notations used in Sect.\  \ref{ssec:Kcomplex}.
We will solve this problem in the following way.
We first observe that Observation 2 in the proof of \cite[Theorem 3.8]{CGRST}
can be interpreted as the following proposition.

\begin{prp}
\label{prp:pr-resol-simp}
For each $a \in \bfL$ and $i \ge 0$,
regard the $i$-th term of the Koszul complex $\calK_a(P_{x})\down$
with $x \in \bfL$ a variable as a covariant functor
$\calK_a(P_{(\blank)})_i \colon \k\bfL \to \mod\k$,
thus as a left $\k\bfL$-module.
Then
the complex $\calK_a(P_{(\blank)})\down$ gives us a
projective resolution of the simple left $\k \bfL$-module corresponding to $a$.
\end{prp}

Denote by $\add \k\bfL$ the formal additive closure of the linear category $\k\bfL$.
By the Yoneda lemma, the covariant Yoneda embedding $\k\bfL \to \mod \k\bfL$ sending
$a \in L$ to $P_a:= \k\bfL(\blank, a)$ induces an equivalence $\add \k\bfL \to \prj \k\bfL$,
by which we sometimes identify these categories.
Moreover, the canonical embedding $\k\bfL \to \add \k\bfL$ induces an equivalence
between their categories of right (resp.\ left) modules.
In particular, this allows us to make the identifications
$\add \k\bfL(\blank, \Ds_{i=1}^m x_i) = \Ds_{i=1}^m \k\bfL(\blank, x_i)$
and
$\add \k\bfL(\Ds_{i=1}^m x_i, \blank) = \Ds_{i=1}^m \k\bfL(x_i, \blank)$
for all $x_1, \dots, x_m \in \bfL$.
We next interpret both the complex $\calK_a(P_{(\blank)})\down$ and the Koszul complex of 
a right $\k\bfL$-module $M$ by using the following notion.

\begin{dfn}
\label{dfn:formalKos_coresol}
Let $a \in \bfL$.
Then the \emph{formal Koszul coresolution} of $a$
is a cocomplex $C\up_a:= (0 \to C_a^0 \ya{d^0} C^1_a \ya{d^1} \cdots)$ in $\add \k\bfL$
defined as follows:
$$
C_a^0:= a,\, C_a^i:= \Ds_{S \in \cov(a)_i} \vee S,\text{ for }i > 0,
$$
and the differential $d^i \colon \Ds_{T \in \cov(a)_{i-1}} \vee T
\to \Ds_{S \in \cov(a)_i} \vee S \ (i > 0)$ is given by the matrix
$$
[\chi(T \subseteq S)(-1)^{n(S,T)} p_{\vee S, \vee T}]_{S, T}
$$
(replace $\vee T$ by $a$ if $i = 1$).
\end{dfn}

Then we immediately see that the following hold.

\begin{prp}
\label{prp:interpret}
The complex $\calK_a(P_{(\blank)})\down$ is nothing but
$(\add \k\bfL)(C_a\up, \blank)$, the image of the contravariant Yoneda imbedding 
of the cocomplex $C_a\up$, and
that the Koszul complex $\calK_a(M)\down$ of $M$ above is obtained by applying
the $M$-dual functor $\Hom_{\k\bfL}(\blank, M)$ to $(\add \k\bfL)(\blank, C_a\up)$, the image
of the covariant Yoneda embedding of $C_a\up$.
\end{prp}

Inspired by this, we define
the (minimal) \emph{$\calI$-relative Koszul coresolution} $\calK\up(V_I)$ of $V_I$
(Definition \ref{dfn:rel-Kos-cores})
as a cocomplex that is sent by the $G$-dual functor $\Hom_A(\blank, G)$ to
the minimal projective resolution of the simple left module corresponding to $V_I$
over the endomorphism algebra of $G$, similar to
$(\add \k\bfL)(\blank, C_a\up)$ above.
This defines the desired (minimal) $\calI$-ralative Koszul complex
$\calK_{I}(M)\down$ of $M$ at $V_I$ (again Definition \ref{dfn:rel-Kos-cores})
by applying the $M$-dual funcctor to $\calK\up(V_I)$ as above.
Thus our $\calK\up(V_I)$ corresponds both to $C\up_a$ in $\add \k\bfL$ and
to $(\add \k\bfL)(\blank, C_a\up)$ in $\prj \k\bfL$.
The proof of our main theorem (Theorem \ref{thm:Kosz-Betti})
verifies the equality \eqref{eq:Koszul-cpx-cdn} by using
the fact that the extension groups
are defined both by a projective resolution of the first variable,
and by an injective resolution of the second variable.

\subsection{Application}
\label{ssec:appl}
In { our} application, we set $A$ to be the incidence algebra of a finite poset $\bfP$,
$\bbI$ to be the set of all intervals in $\bfP$, and $V_I$ to be the interval module
defined by $I$ for all $I \in \bbI$ (Definition \ref{dfn:intv}(2)).
To compute a minimal right/left interval approximation in examples,
we give a handy criterion for a homomorphism to be a right/left interval approximation
(Propositions \ref{prp:crt-r-int-appx}/\ref{prp:crt-r-int-appx}$'$), which
immediately yields a way to construct a minimal right/left interval approximation
(Corollary \ref{cor:form-min-r-appx}/\ref{cor:form-min-r-appx}$'$).
This gives us a way
to check whether a given module $M$ is interval decomposable or not (Definition \ref{dfn:intv}(3)) by using the 0-th interval Betti numbers $\be^0_M(I)$ (Convention \ref{cvn:intv})
only for intervals $I$ such that $V_I$ appears both as a submodule and
a factor module\footnote{%%%%%%%%%%%%%%%%%%%%%%%%%%%%%%%%%%%%%%%%%%%%%%%%%%%
For a submodule $N$ of
an $M \in \mod A$, $M/N$ is sometimes called a \emph{factor} module of $M$.}
of $M$ (Corollary \ref{cor:int-dec-2}),
where the main theorem is used to compute $\be^0_M(I)$.
This problem was also considered in \cite{ABENY} and \cite{AENY-2} that give different ways to check the interval decomposability.
The interval resolutions were also investigated in \cite{BBH}.
On the other hand, when $\bfP$ is a commutative ladder, i.e., a product poset
$\bbA_2 \times \bbA_n\ (n \ge 2)$, where $\bbA_m$ is the totally ordered set with
$m$ elements for each positive integer $m$,
we have a formula that computes the so-called interval replacement of a persistence module
by the interval Betti numbers.
Hence our main theorem makes it possible to compute the interval replacement
by the Koszul coresolution.

\subsection{Plan of the content}
The paper is organized as follows.
In Sect.\ 2, we collect necessary terminologies, such as minimal right (left) $\calI$-approximations.
Especially we give detailed explanation on how to compute minimal right approximations
from right approximations for later use.
In Sect.\ 3, we introduce a concept of relative Koszul coresolution of each indecomposable
direct summand $V_I$ of $G$, which defines a relative Koszul complex of an $A$-module $M$ at $I$,
and prove our main theorem.
In Sect.\ 4 we apply our result to interval resolutions of persistence modules
as explained above, and give some examples.
Finally, in Sect.\ 5, 
we apply our main theorem to the formulas of the ``compressed multiplicities''
and the ``interval replacement'' given in \cite[Theorem 5.5 and Corollary 5.7]{AENY-3}.

\section*{Acknowledgments}
In the original version, $G$ was assumed to be a generator and a cogenerator.
I would like to thank the referee for valuable comments, especially for { their} question
whether the cogenerator assumption on $G$ can be removed.
Thanks to this, I could remove this generator/cogenerator assumption on $G$.

\section{Preliminaries}
{We refer the reader to \cite{ASS} for fundamental facts on
quivers, algebras, and the representation theory of them.}
Throughout this paper, $\k$ is a field, 
$A$ is a finite-dimensional $\k$-algebra,
and all right (rest.\ left) $A$-modules are assumed to be finite-dimensional.
We denote by $D = \Hom_\k(\blank, \k)$ the usual self-duality of $A$,
and by $\mod A$ the category of finite-dimensional right $A$-modules.
In particular, the category of finite-dimensional $\k$-vector spaces is
denoted by $\mod \k$.
The full subcategory of $\mod A$ consisting of the projective (resp.\ injective) modules
is denoted by $\prj A$ (resp.\ $\inj A$).
Moreover, $\{V_I \mid I \in \bbI\}$ is a finite set of mutually non-isomorphic
indecomposable modules in $\mod A$,
and we set $G:= \Ds_{I \in \bbI} V_I$.
Then the endomorphism algebra $\La:= \End_A(G)$ of $G$ turns out to be a \emph{basic} finite-dimensional algebra, namely, as a right $\La$-module, $\La$ is the direct sum of mutually non-isomorphic indecomposable modules.  We regard $G$ as a $\La$-$A$-bimodule.

For a collection $\calS$ of right $A$-modules,
we denote by $\add \calS$ the full subcategory of $\mod A$ consisting of
all direct summands of the finite direct sums of modules in $\calS$.
When $\calS$ consists of only one element, say $\calS = \{M\}$, we simply write
$\add M$ for $\add \{M\}$.
Then in particular, we have $\prj A = \add A$ and $\inj A = \add DA$, where
$A$ and $DA$ are regarded as right $A$-modules.
We set
$\calI:= \add \{V_I \mid I \in \bbI\} = \add G$.

Let $M$ be a right $A$-module.
For a non-negative integer $n$,
the direct sum of $n$ copies of $M$ is usually denoted by $M^n$.
We also use superscript notations for cocomplexes such as
$X = (\cdots \to X^i \ya{d^i} X^{i+1} \to \cdots)$.
We hope that the reader can easily distinguish them from the context.
Recall that $M$ is said to \emph{generate} (resp.\ \emph{cogenerate})
a right $A$-module $N$ if there exists an epimorphism $M^{n} \to N$
(resp.\ a monomorphism $N \to M^{n}$) for some positive integer $n$,
and that $M$ is called a \emph{generator} (resp.\ \emph{cogenerator}) if
it generates (resp.\ cogenerates) all objects in $\mod A$.
As well-known, $M$ becomes a generator (resp.\ cogenerator)
if and only if the right $A$-module
$A$ (resp. $DA$) is contained in $\add M$, which is equivalent to saying that
all indecomposable projective (resp.\ injective) $A$-modules are contained in $\add M$.
In applications, we are interested in the case where $G$ is a generator and a cogenerator,
that is, where $A, DA \in \calI$.

We sometimes use the abbreviation $\dm_A(\blank, ?):= \Hom_A(\blank, ?)$
for variables $\blank, ? $ in $\mod A$.
We denote by $\ell(M)$ the composition length of an $A$-module $M$.
For each morphism $f \colon X \to Y$ in $\mod A$,
we set $\dom(f):= X,\, \cod(f):= Y$.
For a positive integer $n$, we set $[n]:= \{1,\dots, n\}$.
For a set $S$, $\# S$ denotes the cardinality of $S$.

\begin{cvn}
\label{cvn:left-right}
Let $(Q, \ro)$ be a {\em bound quiver}, i.e., a pair of a finite quiver $Q$ and
an admissible ideal $\ro$ of the path algebra $\k Q$ of $Q$,
and suppose that $A$ is defined as $A = \k(Q, \ro):= \k Q/\ro$.
Here, we write the composite of paths in $Q$ \emph{from the left to the right}
(call it \emph{left-to-right notation}),
but that of linear maps between vector spaces \emph{from the right to the left}
(call it \emph{right-to-left notation})
so that representations of the bound quiver $(Q, \ro)$ present {\em right} $A$-modules
as the convention used in \cite{ASS}.

For example the composite of the paths $a_1$ from a vertex 1 to 2 and $a_2$ from 2 to 3
in $Q$ is denoted by $a_1 a_2$, which is a path from 1 to 3.
If $M$ is a representation of $(Q, \ro)$, then the composite of
linear maps $M(1) \ya{M(a_1)} M(2) \ya{M(a_2)} M(3)$
is denoted by $M(a_2) \circ M(a_1)$ defined by
$(M(a_2) \circ M(a_1))(x) = M(a_2)(M(a_1)(x))$
for all $x \in M(1)$, which is equal to $M(a_1 a_2)$.
Therefore, if we use the abbreviation $M(a)(x) = xa$ for a path $a \colon u \to v$ in $Q$
and $x \in M(u)$
(this is used in the proof of Lemma \ref{lem:sub-fac-int}), then the equality
$M(a_2)(M(a_1)(x)) = M(a_1 a_2)(x)$ shows that
$(x a_1)a_2 = x(a_1 a_2)$ for all $x \in M(1)$.
Thus $M$ can be viewed as a right $A$-module, and is naturally extended to
a contravariant functor $\k[Q,\ro] \to \mod \k$,
where $\k[Q,\ro]$ (see Remark \ref{rmk:alg-cat} for details) is a category version of $\k(Q, \ro)$
whose composition is also written by left-to-right notation.

Note that this convention is different from that used in our earlier papers
\cite{ABENY, AENY-2, AENY-3},
where the composite of paths in $Q$ is written by {\em right-to-left} notation.
In particular, with this convention,
the computation of almost split sequences starting from (or ending in)
interval modules are explained in detail in \cite[Sect.\ 5]{ABENY}.
Therefore, care is needed when this is applied in examples in Sect.\ \ref{sec:app-exm}.
\end{cvn}

To { properly} understand the notion of right (resp.\ left) minimal morphisms
that is reviewed below,
we first recall the following well-known category constructions.

\begin{dfn}[Comma categories]
\label{dfn:comm-cats}
Let $M$ be in $\mod A$.
\begin{enumerate}
\item 
We define a category $\calC_M:= \comma{(\mod A)}{M}$ as follows.
The set of objects is given by $(\calC_M)_0:= \bigsqcup_{X \in \mod A}\Hom_A(X, M)$,
and for any $f, f' \in (\calC_M)_0$, we set
$\calC_M(f,f'):= \{g \in \Hom_A(\dom(f), \dom(f')) \mid f = f'g\}$.
The composition of $\calC_M$ is given by that of $\mod A$.

Note that for any $g \in \calC_M(f,f')$, $g$ is an isomorphism in $\calC_M$
if and only if $g$ is an isomorphism in $\mod A$, and that
$\calC_M(f, f') \ne \emptyset$ if and only if $f$ factors through $f'$ in $\mod A$.

\item
We define a category $\calC^M:= \comma{M}{(\mod A)}$ as follows.
The set of objects is given by $\calC^M_0:= \bigsqcup_{X \in \mod A}\Hom_A(M, X)$,
and for any $f, f' \in \calC^M_0$, we set
$\calC^M(f,f'):= \{g \in \Hom_A(\cod(f), \cod(f')) \mid gf = f'\}$.
The composition of $\calC^M$ is given by that of $\mod A$.

Note that for any $g \in \calC^M(f,f')$, $g$ is an isomorphism in $\calC$
if and only if $g$ is an isomorphism in $\mod A$, and that
$\calC^M(f, f') \ne \emptyset$ if and only if $f'$ factors through $f$ in $\mod A$.
\end{enumerate}
\end{dfn}

We cite the following terminologies from \cite{AR1991app} and \cite{ARS}.

\begin{dfn}[Right minimal morphisms]
\label{dfn:r-min}
Let $M$ be in $\mod A$.
\begin{enumerate}
\item
A morphism $f \colon X \to M$ in $\mod A$ is said to be \emph{right minimal}
if for each $h \in \End_A(X)$, the condition $f = fh$ implies that $h$ is an automorphism, or equivalently,
if $\calC_M(f,f)$ consists only of automorpisms of $f$ in $\calC_M$.

\item
Morphisms $f, f' \in \calC_M$ are said to be \emph{equivalent}
if $\calC_M(f,f') \ne \emptyset$ and $\calC_M(f',f) \ne \emptyset$.
The equivalence class containing $f$ is denoted by $[f]$.

\item
Let $f, f' \in \calC_M$.
Then $f'$ is called a \emph{right minimal version} of $f$ if
$f'$ is right minimal and $[f] = [f']$.
\end{enumerate}
\end{dfn}

\begin{rmk}
\label{rmk:r-min}
The following facts are well-known for an $M \in \mod A$.
\begin{enumerate}
\item
By \cite[Proposition 2.1]{ARS},
a morphism $f \colon X \to M$ in $\mod A$ is right minimal if and only if
$\ell(\dom(f))$ is the smallest among $\{\ell(\dom(g)) \mid g \in [f]\}$.
In particular, any $f \in \calC_M$ has its right minimal version,
which is unique up to isomorphism in $\calC_M$.

\item
\cite[Theorem 2.2]{ARS}: For any $f \colon X \to M$ in $\mod A$,
there exists a decomposition $X = X_1 \ds X_2$ such that
$f|_{X_1}$ is right minimal and $f|_{X_2} = 0$.
Moreover, $f|_{X_1}$ is a right minimal version of $f$.

\item
\cite[Corollary 2.3]{ARS}:
A morphism $f \colon X \to M$ is right minimal if and only if
$\Ker f$ contains no nonzero direct summands of $X$.
The latter condition is equivalent to saying that
a section $s \colon X' \to X$ is 0 if $fs = 0$. 

\item
Assume that a morphism $f \colon X \to M$ is decomposed as $f = f'' f'$ in $\mod A$.
Then it is clear by definition that if $f$ is right minimal, then so is $f'$.
\end{enumerate}
\end{rmk}

\begin{dfn}[Right approximations]
\label{dfn:r-appx}
Let $M$ be in $\mod A$.
\begin{enumerate}
\item
A {\em right $\calI$-approximation} of $M$ is a morphism
$f \colon X \to M$ in $\mod A$ with $X \in \calI$ such that
for any $Z \in \calI$, $\dm_A(Z,f) \colon \dm_A(Z, X) \to \dm_A(Z, M)$
is surjective.
Note that if $G$ is a generator, this $f$ necessarily is an epimorphism.
\item
A {\em minimal right $\calI$-approximation} of $M$ is a right $\calI$-approximation of $M$ that is right minimal.
\end{enumerate}
\end{dfn}

\begin{rmk}
\label{rmk:r-appx}
The following facts are well-known for an $M \in \mod A$.
\begin{enumerate}
\item 
Since $\bbI$ is a finite set, there always exist right $\calI$-approximations of $M$.

\item
Let $f, f' \in \calC_M$.
In the case where $\dom(f), \dom(f') \in \calI$ and $\calC(f,f') \ne \emptyset$,
if $f$ is a right $\calI$-approximation of $M$, then so is $f'$. 
In particular, moreover
when $[f] = [f']$, we see that
$f$ is a right $\calI$-approximation of $M$
if and only if so is $f'$.

\item
By (2) above and Remark \ref{rmk:r-min}(2), among right $\calI$-approximations of $M$,
we can always choose a minimal right $\calI$-approximation of $M$, which is
uniquely determined by $M$ up to isomorphism in $\calC_M$ by Remark \ref{rmk:r-min}(1).
\end{enumerate}
\end{rmk}

For a later use in Sect.\ \ref{sec:app-exm},
we here give a way how to compute a minimal right
$\calI$-approximation from a right $\calI$-approximation
(the procedure in the statement (3) below for $\si = \id_{[n]}$),
which would also be helpful to make a computer program.

\begin{lem}
\label{lem:to-r-min}
Assume that a morphism
$f = (f_i)_{i \in [n]} \colon \Ds_{i \in [n]}X_i \to M$ in $\mod A$
is a right $\calI$-approximation with $n$ a positive integer
and $X_i$ indecomposable for all $i \in [n]$, and let $S \subseteq [n]$.
For each $T \subseteq [n]$, we set
$X_T:= \Ds_{i \in T} X_i$ and $f_T:= (f_i)_{i \in T}$ for short.
Then the following are equivalent.
\begin{enumerate}
\item
$f_S \colon X_S \to M$ is a minimal right $\calI$-approximation of $M$.
\item
$f_S$ is a right $\calI$-approximation of $M$, and for each $j \in S$,
$f_{S \setminus \{j\}}$ is not a right $\calI$-approximation of $M$.
\item
There exists some permutation $\si$ of $[n]$ such that
$S = I^\si_n$, where for each $0 \le m \le n$, $I^\si_m$ is defined inductively as follows:
For $m=0$, we set $I^\si_0:= [n]$; and for each $m \in [n]$, we set
\[
I^\si_m:= 
\begin{cases}
I^\si_{m-1} \setminus \{\si(m)\} &
\text{ if $f_{I^\si_{m-1}\setminus \{\si(m)\}}$ is a right
$\calI$-approximation of $M$,}\\
I^\si_{m-1} & \text{ otherwise}.
\end{cases}
\]
\end{enumerate}
In particular, there exists some $S \subseteq [n]$ such that the statement $(1)$ holds.
\end{lem}

\begin{proof}
(3) \implies (2).
Assume that (3) holds, but (2) does not.
Then there exists some $m \in [n]$ with $\si(m) \in I^\si_n$
such that $f_{I^\si_n\setminus \{\si(m)\}}$
is a right $\calI$-approximation of $M$.
Then since $I^\si_n \subseteq I^\si_{m-1}$, we see that
$f_{I^\si_{m-1}\setminus \{\si(m)\}}$ is a right $\calI$-approximation of $M$
by Remark \ref{rmk:r-appx}(2).
Hence $I^\si_n \subseteq I^\si_m = I^\si_{m-1} \setminus \{\si(m)\}$.
Thus $\si(m) \not\in I^\si_n$, a contradiction.

(2) \implies (1).
Assume that (2) holds.
By \cite[Theorem 2.2]{ARS}, there exists a decomposition
$X_S = Y_1 \ds Y_2$ such that the matrix expression of $f_S$ for this decomposition
has the form $f_S = (g_1, 0)$ with $g_1 \colon Y_1 \to M$
a right minimal version of $f$.
Thus $g_1$ is a minimal right approximation of $M$.
It is enough to show that $Y_2 = 0$ because if this is the case,
then $f_S = g_1$, and (1) holds.
Assume contrarily that $Y_2 \ne 0$.
Let $\si_i \colon Y_i \to X_S$ be the inclusions,
and $\pi_i \colon X_S \to Y_i$ the canonical projections for each $i = 1,2$.
Then $f_S = (g_1, 0)$ shows that
\[
g_1 = f_S\si_1,\, f_S\si_2 = 0.
\]
By \cite[12.7 Corollary]{anderson1992rings},
the decomposition $X_S = \Ds_{i \in S}X_i$ of $X_S$ complements $Y_2$,
namely, we have
\begin{equation}
\label{eq:Y2-complement}
X_S = X_T \ds Y_2
\end{equation}
for some $T \subseteq S$.  Since $Y_2 \ne 0$, we have $T \subsetneq S$.
Thus there exists some $j \in S \setminus T$.
The matrix expression of $f_S$ for the decomposition \eqref{eq:Y2-complement}
is given by $f_S = (f_T, 0)$. 
Let $\si' \colon X_T \to X_S$ be the inclusion.
Then $f_T = f_S\si'$ by definition.
Here the direct sum \eqref{eq:Y2-complement} shows that
$\ell:= \pi_1\si' \colon X_T \to Y_1$ is an isomorphism.
Moreover, the diagram
\[
\begin{tikzcd}
\Nname{V'}X_T && \Nname{V}Y_1\\
& \Nname{M}M \\
\Ar{V'}{M}{"f_T" '}
\Ar{V}{M}{"g_1"}
\Ar{V'}{V}{"\ell"}
\end{tikzcd}
\]
is commutative.
Indeed, $f_T - g_1\ell = f_S\si' - f_S\si_1\pi_1\si' = f_S(\id_V - \si_1 \pi_1)\si'
= f_S\si_2 \pi_2 \si' = 0$.
Hence $f_T$ is also a (minimal) right $\calI$-approximation of $M$.
Since $T \subseteq S \setminus \{j\}$, $f_{S \setminus \{j\}}$ is also
a right $\calI$-approximation of $M$, a contradiction.

(1) \implies (3).
Assume that (1) holds.
Let $t$ be the number of elements of $S$.
Then $0 \le t \le n$.
Choose any permutation $\si$ of $[n]$ such that
for each $i \in [n]$,
$\si(i) \in S$ if and only if $n-t < i$.
Then for each $m \in [n - t]$,
we have $I^\si_m = I^\si_{m-1}\setminus \{\si(m)\} \supseteq S$,
and $I^\si_{n-t} = S$.
We show that
$S = I^\si_{m}$ for all $n-t \le m \le n$ by induction on $m$.
Indeed, this is trivial for $m = n-t$.
For $m$ with $n-t < m \le t$, $S = I^\si_{m-1}$ by induction hypothesis.
If $S \ne I^\si_{m}$, then $I^\si_{m} = I^\si_{m-1} \setminus \{\si(m)\} = S \setminus \{\si(m)\}$,
and hence $f_{S \setminus \{\si(m)\}}$ is a right $\calI$-approximation of $M$.
Thus $f_S$ factors through $f_{S \setminus \{\si(m)\}}$.
It is trivial that $f_{S \setminus \{\si(m)\}}$ factors through $f_S$, and hence
we have $[f_S] = [f_{S \setminus \{\si(m)\}}]$.
Then by (1) and Remark \ref{rmk:r-min}(1),
$\ell(X_S) \le \ell(X_{S\setminus \{\si(m)\}})$, a contradiction.
\end{proof}

A \emph{left minimal morphism}, an \emph{equivalence} relation
on the object set of $\calC^M$ with $M \in \mod A$
and a \emph{left minimal version} are defined as dual to those in
Definition \ref{dfn:r-min}.
A \emph{left $\calI$-approximation} and a \emph{minimal left $\calI$-approximation}
are defined as dual versions of Definition \ref{dfn:r-appx}.
The dual versions of Remark \ref{rmk:r-min}, \ref{rmk:r-appx} and
Lemma \ref{lem:to-r-min} also hold, which we omit to state,
but refer to them as Remark \ref{rmk:r-min}$'$, \ref{rmk:r-appx}$'$ and
Lemma \ref{lem:to-r-min}$'$, respectively.

\begin{dfn}
\label{dfn:min-resol-Betti}
Let $M$ be in $\mod A$.
\begin{enumerate}
\item
A sequence
\begin{equation}
\label{eq:min-r-app}
\cdots \ya{f_{r+1}} X_r \ya{f_r} \cdots \ya{f_2} X_1   \ya{f_1} X_0 \ya{f_0} M \ya{f_{-1}} 0
\end{equation}
is called a (resp.\ \emph{minimal}) \emph{$\calI$-resolution} of $M$
if $f_i$ restricts to a (resp.\ minimal) right $\calI$-approximation
$X_i \to \Ker f_{i-1}$
of $\Ker f_{i-1}$ for all $i \ge 0$.
As is easily seen, this sequence is a complex.
Assume the sequence \eqref{eq:min-r-app} is a minimal $\calI$-resolution.
Then it is uniquely determined by $M$ up to isomorphism (of complexes)
by Remark \ref{rmk:r-appx}(3), and
for each $i \ge 0$, $X_i$ has a direct sum decomposition
$X_i \iso \Ds_{I \in \bbI}V_I^{\be^i_M(I)}$ with the unique non-negative integers
$\be^i_M(I)$\ ($I \in \bbI$), which is called the \emph{$\calI$-relative $i$-th Betti number} of $M$ at $I$.

Note that by the definition of right $\calI$-approximations,
the functor $\dm_A(Z, \blank)$ sends the $\calI$-resolution
\eqref{eq:min-r-app} to an exact sequence of vector spaces for all $Z \in \calI$,
and that the sequence \eqref{eq:min-r-app} turns out to be an exact sequnce
if $G$ is a generator.

\item
Dually, a sequence
\begin{equation}
\label{eq:min-l-app}
0 \ya{g^{-1}} M  \ya{g^0} Y^0   \ya{g^1} Y^1\ya{g^2} \cdots \ya{g^r} Y^r \ya{g^{r+1}} \cdots 
\end{equation}
is called a (resp.\ \emph{minimal}) \emph{$\calI$-coresolution} of $M$
if $g^i$ induces a (resp.\ minimal) left $\calI$-approximation
$\Cok g^{i-1} \to Y^i$
of $\Cok g^{i-1}$ for all $i \ge 0$.
As is easily seen, this sequence is a cocomplex.
Assume that the sequence \eqref{eq:min-l-app} is a minimal $\calI$-coresolution.
Then it is uniquely determined by $M$ up to isomorphism (of cocomplexes)
by Remark \ref{rmk:r-appx}$'$(3), and
for each $i \ge 0$, $Y^i$ has a direct sum decomposition
$Y^i \iso \Ds_{I \in \bbI}V_I^{\ovl{\be}^i_M(I)}$ with the unique non-negative integers
$\ovl{\be}^i_M(I)$\ ($I \in \bbI$), which is called the \emph{$\calI$-relative $i$-th co-Betti number} of $M$ at $I$.

Note that by the definition of left $\calI$-approximations,
the functor $\dm_A(\blank, Z)$ sends the sequence
\eqref{eq:min-l-app} to an exact sequence of vector spaces for all $Z \in \calI$,
and that the sequence \eqref{eq:min-l-app} turns out to be
an exact sequence if $G$ is a cogenerator.

\end{enumerate}
\end{dfn}

\section{Relative Koszul coresolutions}

Throughout this section, we fix an $I \in \bbI$, and denote by $e_I \in \La$
the idempotent given by the composite $G \to V_I \to G$ of canonical maps.

\begin{ntn}
We denote the radical of $\mod A$ by $\rad_A$ (\cite[A3.3.3 Definition]{ASS}).
Then $\rad_A(G, G)$ is the Jacobson radical of $\La$ (\cite[A3.3.5(a) Proposition]{ASS}),
and the radical of the left $\La$-module $\dm_A(V_I, G)$ is given by
$\rad_A(V_I, G)$, namely we have
\[
\rad_A(V_I, G) = \rad_A(G,G)\cdot \dm_A(V_I, G)
\]
because the LHS is identified with $\rad_A(G,G)\cdot e_I$ and 
the RHS is identified with $\rad_A(G,G)\cdot \La e_I$.
Similarly, the radical of the right $\La$-module $\dm_A(G, V_I)$ is given by
\[
\rad_A(G, V_I) = \dm_A(G, V_I)\cdot\rad_A(G,G).
\]
We set
$\dm_IS:= \dm_A(V_I, G)/\rad_A(V_I, G)$
(resp.\ $S_I:= \dm_A(G, V_I)/\rad_A(G, V_I)$), which is the simple left (resp.\ right)
$\La$-module corresponding to $I$.
\end{ntn}

From now on, we freely use Auslander--Reiten theory, for which
we refer the reader to \cite[Ch.\ 4]{ASS}.
A left (resp.\ right) almost split morphism that is left (resp.\ right) minimal is called
a \emph{source map} (resp.\ \emph{sink map}) for short
as used in \cite[p.55]{Ringel1099}.

\begin{dfn}
\label{dfn:rel-Kos-cores}
A \emph{minimal $\calI$-relative Koszul coresolution} $\calK\up(V_I)$ 
of $V_I$ is a sequence
\begin{equation}
\label{eq:Kos-cores}
0 \to V_I \ya{d^0} X^1 \ya{d^1} X^2 \ya{d^2} \cdots
\end{equation}
of modules in $\calI$ such that
\begin{enumerate}
\item
$d^0$ is a left minimal version of $\et:= f_0\mu_0$, where $\mu_0 \colon V_I \to E_I$ is a source map from $V_I$ and $f_0 \colon E_I \to X_I$ is a minimal left $\calI$-approximation of $E_I$; and

\item
$d_1 = f_1 \mu_1$, where $\mu_1 \colon X^1 \to \Cok d^0$ is the canonical epimorphism and the sequence
\begin{equation}
\label{eq:min-l-cores-Cok-et}
0 \to \Cok d^0 \ya{f_1} X^2 \ya{d^2} X^3 \ya{d^3} \cdots
\end{equation}
is a minimal left $\calI$-coresolution of $\Cok d^0$.
\end{enumerate}
We denote by $\bar{\calK}\up(V_I)$ the cocomplex obtained by replacing
$d^0$ with $\et$ everywhere in the definition of $\calK\up(V_i)$, and call it
a \emph{pre-minimal} $\calI$-relative Koszul coresolution of $V_I$.
Note that both $\calK\up(V_I)$ and $\bar{\calK}\up(V_I)$ are cocomplexes that are
uniquely determined by $V_I$ up to isomorphism
(of cocomplexes).
When the minimality is not important, we sometimes use $\bar{\calK}\up(V_I)$ instead of $\calK\up(V_I)$.

Let $M$ be in $\mod A$.
Then the complex $\calK_I(M)\down:= \Hom_A(\calK\up(V_I), M)\down$
(resp.\ $\bcalK_I(M)\down:= \Hom_A(\bcalK\up(V_I), M)\down$)
in $\mod \k$
is called a \emph{minimal $\calI$-relative Koszul complex}
(resp.\ \emph{pre-minimal $\calI$-relative Koszul complex})
of $M$ at $I$.
\end{dfn}

\begin{ntn}
Since $d^0$ is a left minimal version of $\et$, $X_I$ has a decomposition
$X_I = X^1 \ds X^{1'}$ such that $\et = \sbmat{d^0\\0}$.
We denote by the projection $X_I \to X^1$ and the  inclusion $X^1 \to X_I$
by $r_I$ and $s_I$, respectively so that we have
\begin{equation}
\label{eq:et-d0}
d^0 = r_I \et, \text{ and }\ \et = s_I d^0.
\end{equation}
We set $C^1:= \Cok d^0$, and for each $i \ge 2$,
decompose $d^i$ as the composite
$d^i = f_i \mu_i$  as in the diagram below, where
$\mu_i \colon X^i \to C^i:= \Cok f_{i-1}$ is the canonical epimorphism,
and $f_i \colon C^i \to X^{i+1}$ is a minimal left $\calI$-approximation of $C^i$:
\[
\begin{tikzcd}[column sep=10pt]
0& V_I && X^1 && X^2 && X^3 &&\cdots\\
&E_I&& X_I & C^1 && C^2 && C^3 &&
\Ar{1-1}{1-2}{}
\Ar{1-2}{1-4}{"d^0"}
\Ar{1-2}{2-4}{"\et" '}
\Ar{1-4}{1-6}{"d^1"}
\Ar{1-6}{1-8}{"d^2"}
\Ar{1-8}{1-10}{"d^3"}
\Ar{1-2}{2-2}{"\mu_0"'}
\Ar{2-2}{2-4}{"f_0"'}
\Ar{2-4}{1-4}{"r_I"}
\Ar{1-4}{2-5}{"\mu_1"}
\Ar{2-5}{1-6}{"f_1"'}
\Ar{1-6}{2-7}{"\mu_2"'}
\Ar{2-7}{1-8}{"f_2"'}
\Ar{1-8}{2-9}{"\mu_3"'}
\Ar{2-9}{1-10}{"f_3"'}
\end{tikzcd}.
\]
\end{ntn}
 
\begin{rmk}
Note that $\et$ does not need to be left minimal although both $\mu_0$ and $f_0$ are left minimal.
An example of $\et$ that is not left minimal is given in Example \ref{exm:CL3}(1).

Assume that $G$ is a cogenerator and that $\et$ is left minimal.
Even in this case,
the sequence \eqref{eq:min-l-cores-Cok-et} is not always exact.
In this case,
the sequence \eqref{eq:min-l-cores-Cok-et} is exact, and hence
the sequence $\calK\up(V_I)$ is exact at $X^i$ for all $i \ge 1$, but
not at $V_I$ in general. 
More precisely, if $V_I$ is non-injective, then $\mu_1$ is a monomorphism
as a left map of the almost split sequence starting from $V_I$, and
$f_1$ is also a monomorphism from the beginning.
Thus the sequence is exact, i.e., the complex $\calK\up(V_I)$ is acyclic.
On the contrary if $V_I$ is injective, then $\mu_1$ is the canonical
epimorphism $V_I \to V_I/\soc V_I$, and hence $\Ker \et = \Ker \mu_1 = \soc V_I$,
and the sequence $\calK\up(V_I)$ is not exact at $V_I$.  
\end{rmk}

\begin{rmk}
\label{rmk:eq-dual}
(1) Recall that the $\La$-$A$-bimodule $G$ defines an adjoint pair $L \adj R$
between $\mod \La$ and $\mod A$, where
$$
L:= \blank \ox_\La G \colon \mod \La \to \mod A, \text{ and } 
R:= \dm_A(G, \blank) \colon \mod A \to \mod \La
$$
with a unit $\et \colon \id_{\mod \La} \To RL$
and a counit $\ep \colon LR \To \id_{\mod A}$.
Since $\ep_G \colon \dm_A(G,G)\ox_\La G \to G$ is an isomorphism,
so is $\ep_M$ for all $M \in \add G = \calI$.
In addition,
since $\et_{\La} \colon \La_{\La} \to \dm_A(G, \La \ox_\La G)$ is an isomorphism,
so is $\et_P$ for all $P \in \add \La_\La = \prj \La$.
Thus this adjoint pair restricts to equivalences that are quasi-inverses to each other:
$$
\begin{tikzcd}
\makebox[1em][r]{$\calI$} & \makebox[1em][l]{$\prj \La$.}
\Ar{1-1}{1-2}{"R", bend left}
\Ar{1-2}{1-1}{"L", bend left}
\end{tikzcd}
$$

(2) Recall also that the $\La$-$A$-bimodule $G$ defines $G$-dual functors
$$
(\blank)^R:= \dm_A(\blank, G) \colon \mod A \to \mod \La\op, \text{ and }
(\blank)^L:= \dm_\La(\blank, G) \colon \mod \La\op \to \mod A
$$
with natural transformations
$\ze \colon \id_{\mod A} \To (\blank)^{RL}$ and
$\ze' \colon \id_{\mod\La} \To (\blank)^{LR}$ defined by evaluations:
$\ze_M \colon M \to \dm_\La(\dm_A(M, G), G),\, x \mapsto (f \mapsto f(x))$
for all $M \in \mod A$ and
$\ze'_N \colon N \to \dm_A(\dm_\La(N, G), G),\, y \mapsto (g \mapsto g(y))$
for all $N \in \mod \La$.
Since $\ze_G \colon G \to \dm_\La(\dm_A(G, G), G)$ is an isomorphism,
so is $\ze_M$ for all $M \in \add G = \calI$.
In addition, since $\ze'_\La \colon \La \to \dm_A(\dm_\La(\La, G), G)$ is an
isomorphism, so is $\ze'_P$ for all $P \in \add \dm_\La\La = \prj \La\op$.
Thus these dual functors restrict to dualities (= contravariant equivalences) that are
quasi-inverses to each other:
$$
\begin{tikzcd}
\makebox[1em][l]{$\calI$} & \makebox[1em][l]{$\prj \La\op$.}
\Ar{1-1}{1-2}{"(\blank)^R", bend left}
\Ar{1-2}{1-1}{"(\blank)^L", bend left}
\end{tikzcd}
$$

(3) Moreover, by applying (2) above to the $\La$-$\La$-bimodule $\La$,
we have the $\La$-dual functors
$$
(\blank)^t:= \dm_\La(\blank, \La_\La) \colon \mod \La \to \mod \La\op, \text{ and }
(\blank)^{t'}:= \dm_{\La\op}(\blank, {}_\La \La) \colon \mod \La\op \to \mod \La,
$$
the restrictions of which on the full subcategories of
finitely generated projective modules turn
out to be dualities that are quasi-inverses to each other:
$$
\begin{tikzcd}
\makebox[1em][r]{$\prj \La$} & \makebox[1em][l]{$\prj \La\op$.}
\Ar{1-1}{1-2}{"(\blank)^t", bend left}
\Ar{1-2}{1-1}{"(\blank)^t", bend left}
\end{tikzcd}
$$

(4) We have an isomorphism $R(V)^t \iso V^R$ that is natural in $V \in \calI$.
In other words, there exists the following diagram that is commutative up to natural isomorphism:
$$
\begin{tikzcd}
\calI & \prj\La\\
\prj\La\op
\Ar{1-1}{1-2}{"R"}
\Ar{1-1}{2-1}{"(\blank)^R"'}
\Ar{1-2}{2-1}{"(\blank)^t"}
\end{tikzcd}
$$
Indeed, this is given as the composite of the isomorphisms
$$
R(V)^t =\dm_\La(R(V), R(G)) \ya{L} \dm_A(LR(V), LR(G)) \ya{\dm_A(\ep_V\inv, \ep_G)} \dm_A(V,G) = V^R
$$
that are natural in $V \in \calI$.
\end{rmk}

Using Remark \ref{rmk:eq-dual}(1) and (2),
we can obtain a way to construct minimal projective resolutions
over $\La$ (resp.\ $\La\op$) as follows:

\begin{lem}
\label{lem:min-l-cores--min-prj-res}
Let $f \colon X \to M$ and $g \colon M \to Y$ be morphisms in $\mod A$
with $X, Y \in \calI$.
Then
\begin{enumerate}
\item
{\rm (i)} $f \colon X \to M$ is right minimal in $\mod A$ if and only if
$$
\dm_A(G, f) \colon \dm_A(G, X) \to \dm_A(G, M)
$$
is left minimal in $\mod \La\op$.

{\rm (ii)}
$f \colon X \to M$ is a minimal right $\calI$-approximation of $M$ in $\mod A$
if and only if
$$
\dm_A(G, f) \colon \dm_A(G, X) \to \dm_A(G, M)
$$
is a projective cover of $\dm_A(G, M)$ in $\mod \La$.

{\rm (iii)} \eqref{eq:min-r-app} is a minimal right $\calI$-resolution of $M$
if and only if the sequence
\[
\footnotesize
\cdots  \ya{}  \dm_A(G, X_r)  \ya{\dm_A(G, f_r)} \cdots \ya{\dm_A(G, f_2)} \dm_A(G, X_1) \ya{\dm_A(G, f_1)} \dm_A(G, X_0)  \ya{\dm_A(G, f_0)}  \dm_A(G,M) \to 0   
\]
is a minimal projective resolution of $\dm_A(G, M)$ in $\mod \La$.
\item
{\rm (i)} $g \colon M \to Y$ is left minimal in $\mod A$ if and only if
$$
\dm_A(g, G) \colon \dm_A(Y, G) \to \dm_A(M,G)
$$
is right minimal in $\mod \La\op$.

{\rm (ii)} $g \colon M \to Y$ is a minimal left $\calI$-approximation of $M$ in $\mod A$
if and only if
$$
\dm_A(g, G) \colon \dm_A(Y, G) \to \dm_A(M,G)
$$
is a projective cover of $\dm_A(M,G)$ in $\mod \La\op$.

{\rm (iii)} \eqref{eq:min-l-app} is a minimal left $\calI$-coresolution of $M$
if and only if the sequence
\[
\footnotesize
\cdots  \ya{}  \dm_A(Y^r,G)  \ya{\dm_A(g^r,G)} \cdots \ya{\dm_A(g^2,G)} \dm_A(Y^1,G) \ya{\dm_A(g^1,G)} \dm_A(Y^0,G)  \ya{\dm_A(g^0,G)}  \dm_A(M,G) \to 0   
\]
is a minimal projective resolution of $\dm_A(M, G)$ in $\mod \La\op$.
\end{enumerate}
\end{lem}

\begin{proof}
We prove only the statement (2) because (1) is proved similarly.
(i) Assume $g$ is left minimal.
To show that $g^R$ is right minimal,
let $h_1 \in \End_{\La\op}(Y^R)$, and assume $g^R h_1 = g^R$.
Then by Remark \ref{rmk:eq-dual}(2), we have
$h_1 = h^R$ for some $h \in \End_A(Y)$, and then
$g^R h^R = g^R$ shows that $hg = g$.
Therefore $h$ is an automorphism by the assumption on $g$, and hence so is $h_1 = h^R$.
Thus $g^R$ is right minimal.

Conversely, assume that $g^R$ is right minimal.
To show that $g$ is left minimal, let $h \in \End_A(Y)$, and assume $hg = g$.
Then $g^R h^R = g^R$, and then $h^R$ is an automorphism.
Hence so is $h$ because $(h^R)^L \iso h$.

(ii) Note that $Y \in \calI$ implies that $Y^R:= \dm_A(Y,G) \in \prj \La\op$.
By definition,
$g \colon M \to Y$ is a left $\calI$-approximation of $M$
if and only if
$g^R:= \dm_A(g, G)$ is an epimorphism from a projective module $Y^R$.
Therefore the second assertion follows.

(iii) This is immediate from (ii).
\end{proof}

Using Remark \ref{rmk:eq-dual}(4), we obtain the following.

\begin{lem}
\label{lem:prj2X-X2inj}
There exists an isomorphism
$$
\dm_\La(X, D(\dm_A(V,G))) \iso D(\dm_\La(\dm_A(G, V),X))
$$
that is natural in $X \in \mod \La$ and in $V \in \calI$.
In particular, for each  $X \in \mod \La$ and each cocomplex $K\up$ in $\calI$,
we have the following isomorphism of cocomplexes:
$$
\dm_\La(X, D(\dm_A(K\up,G))) \iso D(\dm_\La(\dm_A(G, K\up),X)).
$$
\end{lem}

\begin{proof}
As is well-known, we have an isomorphism
$$
\al_{X,P} \colon X \ox_\La P^t \to \dm_\La(P,X),\, x \ox f \mapsto xf(\blank)
$$
that is natural in $X \in \mod \La$, and in $P \in \prj\La$.
This defines an isomorphsim
$$
\be_{X,V}:= \al_{X,R(V)} \colon X \ox_\La (R(V))^t \to \dm_\La(R(V),X).
$$
that is natural in $X \in \mod \La$, and in $V \in \calI$.
By using the usual adjunction of Hom and $\ox$, the isomorphism $\be$
and the isomorphism in Remark \ref{rmk:eq-dual}(4), we obtain an isomorphism
that is natural in $X \in \mod \La$ and in $V \in \calI$ as follows:
$$
\begin{aligned}
\dm_\La(X, D(\dm_A(V,G))) &= \Hom_\La(X, \Hom_\k(V^R, \k))
\iso \Hom_\k(X \ox_\La V^R, \k)\\
&= D(X \ox_\La V^R) \iso D(X \ox_\La R(V)^t) \iso D(\dm_\La(R(V),X))\\
&= D(\dm_\La(\dm_A(G, V),X)).
\end{aligned}
$$
\end{proof}

\begin{rmk}
In the above lemma,
if we ignore the naturality in $X$, the statement is equivalent to saying that
for each  $X \in \mod \La$,
the outer hexagon of the following diagram is commutative up to natural isomorphism:
$$
\begin{tikzcd}
& \prj\La & \mod \k\\
\calI &&& \mod\k.\\
& \prj\La\op & \inj \La
\Ar{2-1}{1-2}{"R"}
\Ar{1-2}{1-3}{"{\dm_\La(\blank, X)}"}
\Ar{1-3}{2-4}{"D"}
\Ar{2-1}{3-2}{"(\blank)^R"'}
\Ar{3-2}{3-3}{"D"}
\Ar{3-3}{2-4}{"{\dm_\La(X,\blank)}"}
\Ar{1-2}{3-2}{"(\blank)^t"}
\Ar{3-2}{1-3}{"X \ox_\La\blank"'}
\end{tikzcd}
$$
The commutativity of the left triangle, of the central triangle, and
of the right quadrangle are given by Remark \ref{rmk:eq-dual}(4), by $\al_{X,\blank}$,
and by Hom-$\ox$ adjunction, respectively.
\end{rmk}

The minimal $\calI$-relative Koszul complex of $M$ gives us a way to compute
the $\calI$-relative Betti numbers of $M$ at $I$ as follows.

\begin{thm}
\label{thm:Kosz-Betti}
Let $M \in \mod A$.
Then the following statements hold.
\begin{enumerate}
\item
$\Hom_A(\calK\up(V_I), G)$ gives a minimal projective resolution
\[
\cdots \ya{} \dm_A(X^2, G) \ya{} \dm_A(X^1, G) \ya{} \dm_A(V_I, G)\to {}_IS \to 0
\]
of the simple left $\La$-module ${}_IS$; and
\item
For each $i \ge 0$,
\[
\dim_\k H_i(\calK_{I}(M)\down) = \be^i_M(I).
\]
\end{enumerate}
\end{thm}

\begin{proof}
(1)
Since ${}_IS = \dm_A(V_I, G)/\rad\dm_A(V_I, G) = V_I^R/\rad V_I^R$, we have an exact sequence
\begin{equation}
\label{eq:prj-cov-simple}
0 \to \rad V_I^R \ya{\si} V_I^R \ya{\pi} {}_IS \to 0
\end{equation}
in $\mod \La\op$,
where $\si$ is the inclusion and $\pi$ is the canonical epimorphism,
and $\pi$ turns out to be a projective cover of ${}_IS$.
Since $(\blank)^R:= \dm_A(\blank, G)$ is left exact, and
the sequence $V_I \ya{d^0} X^1 \ya{\mu_1} C^1 \to 0$ is exact,
the first row in the following commutative diagram is exact:
\[
\begin{tikzcd}
0 & (C^1)^R & (X^1)^R & V_I^R\\
&&  & &\rad V_I^R\\
&& X_I^R&E_I^R
\Ar{1-1}{1-2}{}
\Ar{1-2}{1-3}{"{\mu_1^R}"}
\Ar{1-3}{1-4}{"{(d^0)^R}"}
\Ar{1-3}{3-3}{"{r_I^R}"'}
\Ar{3-3}{3-4}{"{f_0^R}"}
\Ar{3-3}{1-4}{"{\et^R}"}
\Ar{3-4}{1-4}{"{\mu_0^R}"'}
\Ar{2-5}{1-4}{"\si"'}
\Ar{3-4}{2-5}{"\de"'}
\end{tikzcd}.
\]
Since $\mu_0$ is a source map,
we have $\Im \mu_0^R = \rad V_I^R$.
Thus $\mu_0^R$ is written as $\mu_0^R = \si \de$, where
$\de$ is the epimorphism obtained by restricting the codomain of $\mu_0^R$
to its image.
Since $f_0$ is a left $\calI$-approximation, $f_0^R$ is an epimorphism.
Thus $\Im \et^R = \Im \mu_0^R = \rad V_I^R$.
Here we have $\Im \et^R r_I^R = \Im \et^R$.
Indeed, since $\et = s_I d^0$ (see \eqref{eq:et-d0}),
we have $\et^R = (d^0)^R s_I^R$.
Then since $s_I^R$ is an epimorphism, we have
$\Im \et^R = \Im (d^0)^R = \Im \et^R r_I^R$, as desired.
Therefore, $\Im \de f_0^R r_I^R = \rad V_I^R$, which means that
$\de f_0^R r_I^R \colon (X^1)^R \to \rad V_I^R$ is an epimorphism.
Hence the equality
\[
\Im \mu_1^R = \Ker (d^0)^R = \Ker (\si \de f_0^R r_I^R) =  \Ker (\de f_0^R r_I^R)
\]
shows the exactness of the sequence
\begin{equation}
\label{eq:prj-cov-rad}
0 \to (C^1)^R \ya{\mu_1^R} (X^1)^R \ya{\de f_0^R r_I^R} \rad V_I^R \to 0
\end{equation}
in $\mod \La\op$.
Since $d^0 \colon V_I \to X^1$ is left minimal,
$(d^0)^R = \si \de f_0^R r_I^R \colon (X^1)^R \to V_I^R$ is right minimal by
Lemma \ref{lem:min-l-cores--min-prj-res}(2)(i).
Hence so is $\de f_0^R r_I^R$ by Remark \ref{rmk:r-min}(4).
Therefore,  $\de f_0^R r_I^R$ in \eqref{eq:prj-cov-rad} is a projective cover of $\rad V_I^R$.

On the other hand by Lemma \ref{lem:min-l-cores--min-prj-res}(2),
the minimal left $\calI$-coresolution \eqref{eq:min-l-cores-Cok-et} of $C^1$ yiels
a minimal projective resolution
\begin{equation}
\label{eq:prj-res}
\cdots \ya{(d^3)^R} (X^3)^R  \ya{(d^2)^R} (X^2)^R  \ya{(f^1)^R} (C^1)^R \to 0
\end{equation}
of $(C^1)^R$ in $\mod \La\op$.
By combining the projective covers \eqref{eq:prj-cov-simple},
\eqref{eq:prj-cov-rad} and the minimal projective resolution \eqref{eq:prj-res},
we obtain the following
commutative diagram with the first row a minimal projective resolution of ${}_IS$:
\[
\begin{tikzcd}[column sep= 7pt]
\cdots & (X^3)^R &&  (X^2)^R &&  (X^1)^R &&  (V_I)^R & {}_IS & 0.\\
&& (C^2)^R && (C^1)^R && \rad V_I^R
\Ar{1-1}{1-2}{}
\Ar{1-2}{1-4}{"{(d^2)^R}"}
\Ar{1-4}{1-6}{"{(d^1)^R}"}
\Ar{1-6}{1-8}{"(d^0)^R"}
\Ar{1-8}{1-9}{"\pi"}
\Ar{1-9}{1-10}{}
\Ar{1-2}{2-3}{"f_2^R"', twoheadrightarrow}
\Ar{2-3}{1-4}{"\mu_2^R", rightarrowtail}
\Ar{1-4}{2-5}{"f_1^R"', twoheadrightarrow}
\Ar{2-5}{1-6}{"\mu_1^R", rightarrowtail}
\Ar{1-6}{2-7}{"{\de f_0^R r_I^R}"', twoheadrightarrow}
\Ar{2-7}{1-8}{"\si", rightarrowtail}
\end{tikzcd}
\]

(2)
By Lemma \ref{lem:min-l-cores--min-prj-res}(1),
$\dm_A(G, \blank)$ sends a minimal $\calI$-resolution
\eqref{eq:min-r-app} of $M$ with
$X_i \iso \Ds_{J \in \bbI}V_J^{\be^i_M(J)}$ for all $i \ge 0$
to the minimal projective resolution
\begin{equation}
\label{eq:min-prj-res-F}
\cdots \to \dm_A(G,X_1) \ya{\dm_A(G,f_1)} \dm_A(G,X_0) \ya{\dm_A(G,f_0)}
\dm_A(G, M) \to 0
\end{equation}
of $F:= \dm_A(G,M)$ in $\mod \La$.
Then by the minimality of this sequence, the functor $\Hom_\La(\blank, S_I)$ sends the complex $\dm_A(G, X\down)$ to
the complex $\Hom_\La(\dm_A(G, X\down), S_I)$ having the form
\[
0 \to \dm_\La(\dm_A(G,X_0), S_I) \ya{0} \dm_{\La}(\dm_A(G, X_1), S_I) \ya{0} \cdots.
\]
Then noting that
$\Hom_{\La}(\dm_A(G, V_J), S_I) \iso
\begin{cases}
\k & (\text{ if } J = I)\\
0 & (\text{ if } J \ne I)
\end{cases}$
for all $J \in \bbI$, we have
\[
\Ext_\La^i(F, S_I) \iso H^i\Hom_\La(\dm_A(G, X\down), S_I) \iso \k^{\be^i_M(I)}.
\]
Hence we have $\be^i_M(I) = \dim_\k \Ext_\La^i(F, S_I)$.
On the other hand,
by the statement (1), $D(\dm_A(\calK\up(V_I), G))$ turns out to be an injective coresolution of $S_I \iso D(\dm_IS)$.
Hence by Lemma \ref{lem:prj2X-X2inj},
we have isomorphisms
\[
\begin{aligned}
\Ext^i_\La(F, S_I) &\iso H^i \Hom_\La(F, D(\dm_A(\calK\up(V_I), G)))\\
&\iso H^i D(\Hom_\La(\dm_A(G, \calK\up(V_I)), F))\\
&\iso DH_i \Hom_\La(\dm_A(G, \calK\up(V_I)), \dm_A(G, M))\\
&\iso DH_i \Hom_A(\dm_A(G, \calK\up(V_I)) \ox_\La G, M)\\
&\iso DH_i \Hom_A(\calK\up(V_I), M)\\
&= DH_i(\calK_{I}(M)\down).
\end{aligned}
\]
As a consequence, we have
$\be^i_M(I) = \dim_\k DH_i(\calK_{I}(M)\down) = \dim_\k H_i(\calK_{I}(M)\down)$.
\end{proof}

For each additive category $\calC$,
we denote by $\calH^+(\calC)$ (resp.\ $\calH^-(\calC)$) the homotopy category
of cocomplexes bounded below (resp.\ above) in $\calC$.
Since $\calI$ is an additive category,
we can regard the minimal Koszul coresolution $\calK\up(V_I)$ of $V_I$ as an object in
$\calH^+(\calI)$.

\begin{cor}
\label{cor:Kos-res}
Let $M \in \mod A$ and $Y\up \in \calH^+(\calI)$.
Then the following are equivalent:
\begin{enumerate}
\item
$Y\up \iso \calK\up(V_I)$ in $\calH^+(\calI)$.
\item
$(Y\up)^R \iso (\calK\up(V_I))^R$
in $\calH^-(\prj \La\op)$.
\item
$(Y\up)^R:= \dm_A(Y\up, G)$ gives a projective resolution of $\dm_I S$.
\end{enumerate}
If one of them holds, then we have
\[
\dim_\k H_i(\dm_A(Y\up, M)) = \be^i_M(I).
\]
\end{cor}

\begin{proof}
The dualities in Remark \ref{rmk:eq-dual}(2) are extended to
the following dualities that are quasi-inverses to each other:
$$
\begin{tikzcd}[column sep=50pt]
\makebox[1em][r]{$\calH^+(\calI)$} & \makebox[1em][l]{$\calH^-(\prj \La\op)$.}
\Ar{1-1}{1-2}{"(\blank)^R", bend left}
\Ar{1-2}{1-1}{"(\blank)^L", bend left}
\end{tikzcd}
$$
Hence the statement (1) is equivalent to the statement (2).
The statement (2) is equivalent to the statement (3) by
Theorem \ref{thm:Kosz-Betti}(1).

Note that the contravariant functor $\dm_A(\blank, M) \colon \calI \to \mod \k$ is extended to
a contravariant functor $\calH^+(\calI) \to \calH^-(\mod\k)$.
Then the statement (1) implies that
$\dm_A(Y\up, M)\iso \dm_A(\calK\up(V_I), M) = \calK_I(M)\down$
in $\calH^-(\mod \k)$, which shows that
$H_i(\dm_A(Y\up, M)) \iso H_i(\calK_I(M)\down)$ in $\mod \k$ for all $i \ge 0$.
Hence the remaining assertion follows by Theorem \ref{thm:Kosz-Betti}(2).
\end{proof}

\begin{dfn}
A cocomplex $Y\up \in \calH^+(\calI)$ is called
an \emph{$\calI$-relative Koszul coresolution} of $V_I$
if $Y\up$ has the form
$Y\up = (0 \to V_I \ya{d^0_Y} Y^1 \ya{d^1_Y} \cdots)$, and
satisfies one of the equivalent conditions in Corollary \ref{cor:Kos-res}.
If this is the case, then for each $M \in \mod A$, we call $\dm_A(Y\up, M)$
an $\calI$-relative \emph{Koszul complex} of $M$ at $I$.
\end{dfn}

\begin{rmk}
Using the condition (3) in Corollary \ref{cor:Kos-res}, we see that the following
are equivalent:
\begin{enumerate}
\item 
$Y\up$ is an $\calI$-relative Koszul coresolution of $V_I$.
\item
$d^0_Y$ is equivalent to $g_0\nu_0$ in the comma category
$\calC^{V_I}$,
where $\nu_0 \colon V_I \to E_{Y,I}$ is a left almost split morphism from $V_I$,
$g_0$ is a left $\calI$-approximation of $E_{Y,I}$, and
$d^1_Y = g_1\nu_1$, where
$\nu_1$ is a cokernel morphism of $d^0_Y$, and the sequence
\begin{equation}
\label{eq:l-cores-Cok-et}
0 \to \Cok d^0_Y \ya{g_1} Y^2 \ya{d^2_Y} Y^3 \ya{d^3_Y} \cdots
\end{equation}
is a left $\calI$-coresolution of $\Cok d^0_Y$.
\end{enumerate}
The statement (2) above is equivalent to saying that
we have a commutative diagram of the following form:
\[
\begin{tikzcd}[column sep=10pt]
0& V_I && Y^1 && Y^2 && Y^3 &&\cdots,\\
&E_{Y,I}&& Y_I & \Cok d^0_Y && \Cok d^1_Y && \Cok d^2_Y &&
\Ar{1-1}{1-2}{}
\Ar{1-2}{1-4}{"d^0_Y"}
\Ar{1-4}{1-6}{"d^1_Y"}
\Ar{1-6}{1-8}{"d^2_Y"}
\Ar{1-8}{1-10}{"d^3_Y"}
\Ar{1-2}{2-2}{"\nu_0"'}
\Ar{2-2}{2-4}{"g_0"'}
\Ar{2-4}{1-4}{"r_Y", shift left=3pt}
\Ar{1-4}{2-4}{"s_Y", shift left=3pt, dashed}
\Ar{1-4}{2-5}{"\nu_1"}
\Ar{2-5}{1-6}{"g_1"'}
\Ar{1-6}{2-7}{"\nu_2"'}
\Ar{2-7}{1-8}{"g_2"'}
\Ar{1-8}{2-9}{"\nu_3"'}
\Ar{2-9}{1-10}{"g_3"'}
\end{tikzcd}.
\]
where for each $i \ge 1$,
$\nu_i$ is the cokernel morphism of $d^{i-1}_Y$,
$g_i$ is a left $\calI$-approximation of $\Cok d^{i-1}_Y$.
The commutativity of the left square means that it is commutative with $r_Y$
and $s_Y$, respectively.

Hence for instance, the minimal (resp.\ pre-minimal) $\calI$-relative Koszul coresolution is
an $\calI$-relative Koszul coresolution.

By Corollary \ref{cor:Kos-res}, we see that
to compute the relative Betti numbers, the minimality of the Koszul coresolution
is not necessary, just an $\calI$-relative Koszul coresolution is enough for this purpose.
\end{rmk}

Consider the case where
$\La$ is isomorphic to the incidence algebra of a finite lower semi-lattice $\bfL$,
or equivalently, the full subcategory $\calP$ of $\mod A$ consisting of the set
$\{V_I \mid I \in \bbI\}$ of objects is isomorphic to the incidence category $\k\bfL$,
and identify $\calP$ with $\k\bfL$.
Then we can come back to the original situation, and
Corollary \ref{cor:Kos-res} gives us a way to compute a Koszul coresolution
as in the original way explained in the introduction.
For instance, this occurs in the following case.

\begin{exm}
\label{exm:lattice-An}
For Auslander--Reiten quivers ({\em AR-quivers} for short),
we refer the reader to \cite[Ch.\ 4]{ASS} (see also \cite[2.3]{Ringel1099}),
but we treat them as translation quivers as in \cite[p.61]{Ringel1099}.

Let $A$ be the incidence algebra of the poset of type $\bbA_n$
for some $n \ge 2$ (see Sect.\ \ref{ssec:appl}), and $\calP = \{V_I \mid I \in \bbI\}$
be given as
$\{V \in  \Ga_A \mid \Hom_A(M, V) \ne 0\}$ for some $M \in \Ga_A$, where $\Ga_A$ is
the AR-quiver of $A$.
Then $\calP$ is isomorphic to the incidence category of a finite lattice $\bfL$
(in this case, it becomes a 2D-grid),
where the full subquiver of $\Ga_A$ consisting of modules in $\calP$
gives the Hasse quiver of $\bfL$.
For example, if $n = 4,\, \bbA_4 = (1 \to 2 \to 3 \to 4)$ and $M$ has dimension vector
$(0,0,1,1)$, then $\Ga_A$ and $\calP$ is displayed as follows:
\[\begin{tikzcd}
	&&& \Nname{top}\circ \\
	&& \circ && \Nname{right}\circ \\
	&\Nname{M} M && \circ && \circ \\
	\circ && \Nname{bottom}\circ && \circ && \circ
	\arrow[from=1-4, to=2-5]
	\arrow[from=2-3, to=1-4]
	\arrow[dashed, no head, from=2-3, to=2-5]
	\arrow[from=2-3, to=3-4]
	\arrow[from=2-5, to=3-6]
	\arrow[from=3-2, to=2-3]
	\arrow[dashed, no head, from=3-2, to=3-4]
	\arrow[from=3-2, to=4-3]
	\arrow[from=3-4, to=2-5]
	\arrow[dashed, no head, from=3-4, to=3-6]
	\arrow[from=3-4, to=4-5]
	\arrow[from=3-6, to=4-7]
	\arrow[from=4-1, to=3-2]
	\arrow[dashed, no head, from=4-1, to=4-3]
	\arrow[from=4-3, to=3-4]
	\arrow[dashed, no head, from=4-3, to=4-5]
	\arrow[from=4-5, to=3-6]
	\arrow[dashed, no head, from=4-5, to=4-7]
\ar[to path={(M.west)--(top.north)--
(right.east)--(bottom.south)--(M.west)}, dash]
\end{tikzcd},
\]
where the quadrangle stands for $\calP$.

Another example is given in Example \ref{exm:La=poset}.
\end{exm}

\begin{cor}
\label{cor:weakKos_cores}
Assume the setting as above, and 
identify $\calI = \add \calP$ with the formal additive closure $\add \k\bfL$ of $\k\bfL$.
Then the formal Koszul coresolution
$C\up_{I}$ in $\calI$  {\em (see Definition \ref{dfn:formalKos_coresol})}
is an $\calI$-relative Koszul coresolution,
and hence for each $i \ge 0$, we have
$$
\dim H_i(\dm_A(C_{I}\up, M)) = \be_M^i(I).
$$
\end{cor}

\begin{proof}
By Propositions \ref{prp:pr-resol-simp} and  \ref{prp:interpret},
we see that $\dm_A(C\up_{I}, G)$ gives a projective resolution of $\dm_I S$.
Hence Corollary \ref{cor:Kos-res} proves the assertion.
\end{proof}

The following is immediate from Theorem \ref{thm:Kosz-Betti}.

\begin{cor}
\label{cor:int-dec-chrctrz}
The following are equivalent for all $M \in \mod A$:
\begin{enumerate}
\item
$M \in \calI$;
\item
$H_1(\calK_{I}(M)\down) = 0$ for all $I \in \bbI$;
\item
$\dim_\k M = \sum_{I \in \bbI} \dim_\k H_0(\calK_{I}(M)\down)\dim_\k V_I$; and
\item
$M \iso \Ds_{I \in \bbI}V_I^{\dim_\k H_0(\calK_{I}(M)\down)}$.
\end{enumerate}
\end{cor}

\begin{proof}
Each statement is equivalent to the fact that $f_0$ in the minimal $\calI$-resolution
\eqref{eq:min-r-app} is an isomorphism.
\end{proof}

The following is useful in applications.

\begin{rmk}
To compute $H_0(\calK_{I}(M)\down) = H_0(\bcalK_{I}(M)\down)$
in (3) and (4) above, we only need to know
the beginning part (call it the {\em 0-1-part})
\[
0 \to V_I \ya{\et} X^1
\quad\text{and}\quad
\dm_A(X^1, M) \ya{\dm_A(\et, M)} \dm_A(V_I, M) \to 0
\]
of $\bcalK\up(V_I)$ with $\et = f_0\mu_0$ and of $\bcalK_{I}(M)\down$,
which gives us the formula
\begin{equation}
\label{eq:beta0-formula}
\be^0_M(I) = \dim \Cok \dm_A(\et,M) = \dim \dm_A(V_I, M) - \dim \Im \dm_A(f_0\mu_0, M).
\end{equation}
Therefore, to apply Corollary \ref{cor:int-dec-chrctrz} in the setting of Sect.\ 4,
it would be convenient first to compute the 0-1-parts of Koszul coresolutions
for all interval modules $V_I\ (I \in \bbI)$ (Definition \ref{dfn:intv}).
(Note that we can further restrict $I \in \bbI$ as in Corollary \ref{cor:int-dec-2}
and in Proposition \ref{prp:smaller-Sint}.)
To this end, we need to compute the source map $\mu_0 \colon V_I \to E_I$
and the minimal left interval approximation $f_0 \colon E_I \to X^1$ for all $I \in \bbI$.
These can be done by using a recipe explained in \cite[Sect.\ 5]{ABENY}
and Corollary \ref{cor:form-min-r-appx}$'$ in the next section, respectively.
(We may replace $f_0$ by just a left interval approximation of $E_I$, and in this case
we can apply Corollary \ref{cor:const-r-int-appx}'.)
\end{rmk}

\section{Applications and examples}
\label{sec:app-exm}

Throughout this section, $Q$ is an acyclic quiver without multiple arrows,
$\ro$ is the ideal of the path-algebra $\k Q$ generated by the full commutativity relations, and set $A:= \k Q/\ro$.
Then the vertex set $Q_0$ turns out to be a poset $(Q_0, \preceq)$ by the partial order $\preceq$ defined by $x \preceq y$ if and only if there exists a path from $x$ to $y$ in $Q$ for all $x, y \in Q_0$.

\begin{dfn}
\label{dfn:incidence}
Let $S = (S, \le)$ be a finite poset.
\begin{enumerate}
\item
For any $x, y \in S$ with $x \le y$, we set $[x,y]:= \{i \in S \mid x \le i \le y\}$
and call it the {\em segment} from $x$ to $y$.
The set $\{[x,y] \mid x \le y, x, y \in S\}$ of all segments is denoted by $\Seg(S)$.
\item
The {\em incidence algebra} $\k S$ of $S$ over $\k$
is defined to be the algebra having $\Seg(S)$ as a basis
with the multiplication defined by $[x,y][u,v]:= \de_{y,u}[x,v]$ for all
$[x,y] , [u,v] \in \Seg(S)$.
\item
The {\em Hasse quiver} $H(S)$ of $S$ is defined as follows:
$H(S)_0:= S$, $H(S)_1:= \{a_{x,y} \mid x, y \in S, x < y, [x,y] = \{x, y\}\}$, and
the source and the target of $a_{x,y}$ is $x$ and $y$,
respectively for all $a_{x,y} \in H(S)_1$.
Then $H(S)$ is an acyclic quiver without multiple arrows.
The {\em Hasse bound quiver} $(H(S), \ro)$ is the pair of $H(S)$ and the ideal $\ro$
of the path algebra $\k H(S)$ generated by the full commutativity relations.
\end{enumerate}
\end{dfn}

\begin{rmk}
\label{rmk:incidence}
Let $S = (S, \le)$ be a finite poset, and
$(H(S), \ro)$ the Hasse bound quiver of $S$.
\begin{enumerate}
\item
There exists a bijection from the set of paths in $H(S)$ modulo $\ro$ to
the set $\Seg(S)$ preserving the source and the target that induces an algebra isomorphism
$\k(H(S), \ro) \to \k S$.
\item
In particular, for the poset $S:= (Q_0, \preceq)$,
note that $(Q, \ro)$ is isomorphic to $(H(S), \ro)$, and hence $A$ is isomorphic to
the incidence algebra $\k S$.
\end{enumerate}
\end{rmk}

\begin{dfn}[intervals, interval modules, interval decomposable modules]
\label{dfn:intv}
\hspace{1pt}
\begin{enumerate}
\item
A full subquiver $I$ of $Q$ is called an \emph{interval} if it is connected and \emph{convex} in the sense that if $x, y$ are vertices in $I_0$ and $p$ is a path in $Q$ from $x$ to $y$, then $p$ is a path in $I$.

We denote by $\bbI$ the set of all interval subquivers of $Q$.
\item
Each interval subquiver $I \in \bbI$ defines an indecomposable representation $V_I$ of $Q$ by
setting $V_I(x) = \k$ (resp.\ $V_I(x) = 0$) if $x \in I_0$ (resp.\ $x \not\in I_0$)
for all $x \in Q_0$, and
$V_I(a) = \id_\k$ (resp.\ $V_I(a) = 0$) if $a \in I_1$ (resp.\ $a \not\in I_1$)
for all $a \in Q_1$.

For each $I \in \bbI$, note that $V_I$ satisfies all the commutativity relations on $Q$,
and hence can be seen as a right $A$-module.
A right $A$-module is called an \emph{interval module}
if it is isomorphic to $V_I$ for some $I \in \bbI$.
\item
A right $A$-module is said to be
\emph{interval decomposable} if it is the direct sum of some
interval modules.

We set $G:= \Ds_{I \in \bbI}V_I$ and $\calI:= \add G$.
Then a right $A$-module is interval decomposable if and only if it is in $\calI$.
\end{enumerate}
\end{dfn}

Note that indecomposable projective modules and indecomposable injective modules
are interval modules, and hence $G$ is a generator and a cogenerator in $\mod A$.
Therefore, we can apply all the results in Section 3 in this setting.
Note also in this case that $\La$ has finite global dimension by \cite[Proposition 4.5]{AENY-3},
and hence all minimal interval
resolutions (resp.\ coresolutions) are bounded complexes (resp.\ cocomplexes).

\begin{cvn}
\label{cvn:intv}
Here, ``$\calI$-relative'' or ``$\calI$-'' is replaced with the word ``interval''.
For instance, a minimal $\calI$-resolution and an $\calI$-relative Koszul coresolution are called a \emph{minimal interval resolution} and an \emph{interval Koszul coresolution}, respectively.
\end{cvn}

To make the computations of minimal right (resp.\ left) interval approximations easier,
we first give a handy criterion to check whether a homomorphism is a right (resp.\ left) interval approximation.
By using this, we give examples of the minimal interval Koszul coresolution of an interval module
$V_I$, the minimal interval Koszul complex $\calK_{I}(M)\down$ of a module $M$, and
compute some interval Betti numbers $\be^i_M(I)$.

\subsection{A criterion for right/left interval approximation}

We start with an easy remark.

\begin{rmk}
\label{rmk:r-int-appx}
Let $f \colon X \to M$ be in $\mod A$.
It is clear from Definition \ref{dfn:r-appx}
that the following are equivalent:
\begin{enumerate}
\item
$f$ is a right interval approximation of $M$.
\item
For any $I \in \bbI$ and any $g \in \Hom_A(V_I, M)$,
$g = fh$ for some $h \colon V_I \to X$.
\end{enumerate}
\end{rmk}

We want to make the sets $\bbI$ and $\Hom_A(V_I, M)$ in the statement (2) above
much smaller.

\begin{ntn}
For any $X, Y \in \mod A$,
we denote by $\Mon(X,Y)$ (resp.\ $\Epi(X, Y)$)
the set of all monomorphisms (resp.\ epimorphisms) in $\Hom_A(X, Y)$.
\end{ntn}

\begin{dfn}
Let $M$ be in $\mod A$.
Then we set
\[
\begin{aligned}
\Sint(M)&:= \{I \in \bbI \mid \Mon(V_I, M) \ne \emptyset\},\\
\Fint(M)&:= \{I \in \bbI \mid \Epi(M, V_I) \ne \emptyset\}.
\end{aligned}
\]
\end{dfn}

The following is a special property of interval decomposable modules,
the first half of which was used to show the global dimension of $\La$ is finite in \cite{AENY-3}.
The second half is used in the proof of Proposition \ref{prp:crt-r-int-appx}.

\begin{lem}
\label{lem:sub-fac-int}
Let $I \in \bbI$.
Then all submodules and all factor modules of $V_I$ are interval decomposable.
\end{lem}

\begin{proof}
For the submodules this statement was already proved in \cite{AENY-3}.
The statement for factor modules is proved similarly as follows.
Let $W$ be a factor module of $V_I$.
Then $W = V_I/V'$ for some $V' \le V_I$.
Denote $1 \in \k = V_I(x)$ by $1_x$ for all $x \in Q_0$.

We show that the full subquiver $\supp W$ of $Q$ with the vertex set
$(\supp W)_0:= \{x \in Q_0 \mid W(x) \ne 0\}$ is convex.
Let $x, y \in (\supp W)_0$, and $p$ a path from $x$ to $y$ in $Q$,
and suppose that $p = p'p''$ for some paths $p'$ from $x$ to a vertex $z$,
and $p''$ from $z$ to $y$.
To show the assertion, it is enough to show that $z \in (\supp W)_0$.
Assume to the contrary that $z \not\in (\supp W)_0$.
Then $1_z \in V'(z)$, and $1_y = 1_zp'' \in V'(y)$, and hence $W(y) =0$, a contradiction.
Thus $\supp W$ must be convex.
Let $C_1, \dots, C_t$ be the connected components of $\supp W$.
Then they are connected and convex, and $W = W|_{C_1} \ds \cdots \ds  W|_{C_t}$,
where $W|_{C_i}$ are the restrictions of $W$ on $C_i$ for all $i = 1,\dots, t$.
As is easily seen, $W|_{C_i}$ are interval modules for all $i$.
Hence $W$ is interval decomposable.
\end{proof}

Remark \ref{rmk:r-int-appx} and Lemma \ref{lem:sub-fac-int} immediately proves the following.

\begin{prp}
\label{prp:crt-r-int-appx}
Let $f \colon X \to M$ be in $\mod A$.
Then the following are equivalent.
\begin{enumerate}
\item
$f$ is a right interval approximation of $M$.
\item
For any $I \in \Sint(M)$ and any monomorphism $g \colon V_I \to M$ in $\mod A$,
$g = fh$ for some $h \colon V_I \to X$.
\end{enumerate}
\end{prp}

\begin{proof}
(1) \implies (2). This is trivial.

(2) \implies (1). Let $I \in \bbI$ and $g \colon V_I \to M$ be in $\mod A$.
Since $\Im g$ is a factor module of $V_I$, it is interval decomposable
by Lemma \ref{lem:sub-fac-int}, thus we have an isomorphism
$\Im g \iso \Ds_{i=1}^n V_{I_i}$
for some $I_i \in \bbI$\ ($i = 1,\dots, n$), by which we identify these modules.
Consider the factorization $g = \si g'$ of $g$ by the restriction $g' \colon V_I \to \Im g$
of $g$ and the inclusion $\si \colon \Im f \to M$.
Let $g' = \dm^t(g'_i)_{i=1}^n$ and $\si = (\si_i)_{i=1}^n$ be the matrix presentations
of $g'$ and $\si$ with this decomposition of $\Im g$, respectively.
Then $\si_i \colon V_{I_i} \to M$ is a monomorphism with $I_i \in \Sint(M)$
for all $i = 1,\dots, n$.
Therefore by (2), for each $i$, there exists some $h_i \in \dm_A(V_{I_i}, X)$ such that
 $\si_i = f h_i$.
 Hence $g = \sum_{i=1}^n \si_i g'_i = f(\sum_{i=1}^n h_i g'_i)$.
\end{proof}

Proposition \ref{prp:crt-r-int-appx} above immediately gives us a way to
construct a right interval approximation as follows.

\begin{cor}
\label{cor:const-r-int-appx}
Let $M \in \mod A$.
For each $I \in \Sint(M)$, let $W_I$ be the subspace of the finite-dimensional
vector space $\Hom_A(V_I, M)$ spanned by $\Mon(V_I, M)$.
Then $n_I:= \dim W_I$ is finite, and hence
there exists a finite subset $\{\bar{f}_I^{(1)}, \dots, \bar{f}_I^{(n_I)}\}$ of $\Mon(V_I,M)$
that is a basis of $W_I$, thus $W_I = \Ds_{i=1}^{n_I} \k \bar{f}_I^{(i)}$.
Then
\[
\bar{f}:= ((\bar{f}_I^{(1)}, \dots, \bar{f}_I^{(n_I)}))_{I \in \Sint(M)} \colon \Ds_{I \in \Sint(M)} V_I^{n_I} \to M
\]
is a right interval approximation of $M$
\end{cor}

\begin{proof}
Let $I \in \Sint(M)$ and $g \in \Mon(V_I, M)$.
Then by construction, we have $g = \sum_{i=1}^{n_I} k_i \bar{f}_I^{(i)}$
for some $k_i \in \k\ (i = 1, \dots, n_I)$.
By taking
\[
h:= {}^t({}^t(h_J^{(1)}, \dots, h_J^{(n_J)}))_{J \in \Sint(M)} \in \Hom_A(V_I, \Ds_{I \in \Sint(M)} V_I^{n_I}),
\]
where $h_I^{(i)}:= \de_{J,I} k_i\id_{V_I}$, we have $g = \bar{f}h$.
Therefore, $\bar{f}$ is a right interval approximation.
\end{proof}

Corollary \ref{cor:const-r-int-appx} and Lemma \ref{lem:to-r-min} gives us
a way to construct a \emph{minimal} right interval approximation as follows.

\begin{cor}
\label{cor:form-min-r-appx}
Let $M \in \mod A$, and construct a right interval approximation $\bar{f}$ of $M$
as in {\rm Corollary \ref{cor:const-r-int-appx}}.
Choose a finite subset $S \subseteq \Sint(M)$ and 
a subset $\{f_I^{(1)}, \dots, f_I^{(m_I)}\}$
of $\{\bar{f}_I^{(1)}, \dots, \bar{f}_I^{(n_I)}\}$ for each $I \in S$
by the procedure described in
{\em Lemma \ref{lem:to-r-min}(3)} applied to the identity permutation $\si$.
Then $\{f_I^{(i)}\}_{i \in [m_I]}$
is a linearly independent subset of $\Mon(V_I, M)$ such that
\[
f:= ((f_I^{(1)}, \dots, f_I^{(m_I)}))_{I \in S} \colon \Ds_{I \in S} V_I^{m_I} \to M
\]
is a minimal right interval approximation of $M$.
\qed
\end{cor}

\begin{rmk}
After publishing the preprint version of this paper to arXiv,
the preprint of the paper \cite{aoki2023summand} appeared,
in which Theorem 3.4 states the following:
Let $A$ be a finite-dimensional $\k$-algebra, $\calX$ a full subcategory of $\mod A$
with $\calX = \add \calX$ that contains any quotients of
the indecomposables in $\calX$,
and $f \colon X \to M$ a minimal right $\calX$-approximation.
Then 
\begin{enumerate}
\item
for any indecomposable direct summand $X'$ of $X$, the restriction
$f|_{X'}$ is injective.
\item
$\supp X \subseteq \supp M$.
\item
If $\calX$ contains $A$, then $\supp X = \supp M$.
\end{enumerate}
The statement (1) is the essential part, and (2), (3) are immediate from (1) and Definition \ref{dfn:r-appx}.
Note that (2), (3) also follows from the existence of an $f$
in Corollary \ref{cor:form-min-r-appx}; the existence of this $f$ 
follows from (1)
(except for the linear independence of $\{f_I^{(i)}\}_{i \in [m_I]}$, which follows
from Lemma \ref{lem:to-r-min}).
However, this theorem (even its interval version \cite[Corollary 3.11]{aoki2023summand}) gives only necessary conditions for a morphism $f$ to be a minimal right $\calX$-approximation of $M$
and does not give a way to construct it.
Conversely, note that also in this setting
(at least when $\calX$ contains only finitely many
indecomposables up to isomorphism),
Corollary \ref{cor:form-min-r-appx} holds by the same argument used in its proof,
and it gives a construction of a minimal right $\calX$-approximation.

We here note that 
Remark \ref{rmk:r-min}(1) gives us an alternative proof of (1) as follows:
By assumption, we have a decomposition $X = X' \ds X''$ with $X', X'' \in \calX$.
Set $(f',f''):= (f|_{X'}, f|_{X''}) = f$, and let
$f'_2 \colon f(X') \hookrightarrow M$ be the inclusion.
Then 
$(f'_2, f'') \colon f(X') \ds X'' \to M$ is equivalent to $f$.
(Indeed,
since $f' = f'_2 f'_1$ with the restriction $f'_1 \colon X' \to f(X')$ of $f'$,
we have $f = (f'_2, f'') \binom{f'_1}{\id_{X''}}$.
By assumption, we have $f(X')\ds X'' \in \calX$, and
hence $(f'_2, f'')$ factors through $f$ because $f$ is a right $\calX$-approximation of $M$.
Thus $[f] = [(f'_2, f'')]$ as desired.)
Then by Remark \ref{rmk:r-min}(1), we have 
$\ell(X'\ds X'') =\ell(X) \le \ell(f(X') \ds X'')$,
and hence $\ell(X') \le \ell(f(X'))$. 
Thus $f|_{X'}$ is injective.

Both Proposition \ref{prp:crt-r-int-appx} and \cite[Theorem 3.4]{aoki2023summand} use the fact that $\calX$ is closed under
taking factor modules of indecomposables in $\calX$
in the same way, but
there is a difference between them:
the former goes to a necessary and sufficient condition for a morphism 
$f$ to be a right $\calX$-approximation, and the latter goes to 
a necessary condition for $f$ to be a minimal right $\calX$-approximation.

Nevertheless, the condition (3) above
can be used as a handy criterion to check whether
$f_{I^\si_{m-1} \setminus \{\si(m)\}}$ is not a right $\calX$-approximation
in Lemma \ref{lem:to-r-min}(3).
\end{rmk}

We can further make the set of necessary intervals $I$ smaller than $\Sint(M)$.
To this end, we introduce a partial order on $\Sint(M)$ that may different from the
inclusion order.

\begin{dfn}
For a poset $S$ and a subset $T$ of $S$, $T$ is said to be an \emph{up-set} in $S$
if for any $a,\, b \in S$ with $a \le b$, $a \in T$ implies $b \in T$.
Dually, $T$ is said to be a \emph{down-set} in $S$ if for any
$a,\, b \in S$ with $a \le b$, $b \in T$ implies $a \in T$.

Note that for any $I,\, J \in \bbI$, there exists a monomorphism $V_I \to V_J$
(resp.\ an epimorphism $V_J \to V_I$) if and only if
$I$ is an up-set (resp.\ down-set) in $J$.
When $I$ is an up-set in $J$, we choose the canonical monomorphism
$\si:= \si_{J,I} \colon V_I \to V_J$ defined by setting
$\si_x \colon V_I(x) \to V_J(x)$ to be the identity map of $\k$ if $x \in I$,
and 0 otherwise, among monomorphisms from $V_I$ to $V_J$.
Dually, when $I$ is a down-set in $J$, we choose the canonical epimorphism
$\pi:= \pi_{I,J} \colon V_J \to V_I$ defined by setting
$\pi_x \colon V_J(x) \to V_I(x)$ to be the identity map of $\k$ if $x \in I$,
and 0 otherwise, among epimorphisms from $V_J$ to $V_I$.
\end{dfn}

\begin{dfn}
\label{dfn:order-S}
Let $M \in \mod A$. Then we define a partial order $\le_S$ on $\Sint(M)$ as follows:
For any $I,\, J \in \Sint(M)$,
$I \le_S J$ if and only if $I$ is an up-set in $J$,
and any monomorphism $V_I \rightarrowtail M$ is extended to a monomorphism $V_J \rightarrowtail M$
along $\si_{J,I}$ (see the diagram \eqref{eq:Sint-ext}).
We set
\[
\max \Sint(M):= \{I \in \Sint(M) \mid I \text{ is maximal in } (\Sint(M), \le_S)\}.
\]
\end{dfn}

Using the partial order $\le_S$ in $\Sint(M)$, Proposition \ref{prp:crt-r-int-appx}
is improved as follows:

\begin{prp}
\label{prp:smaller-Sint}
Let $f \colon X \to M$ be in $\mod A$.
Then the following are equivalent.
\begin{enumerate}
\item
$f$ is a right interval approximation of $M$.
\item[(2$'$)]
For any $J \in \max\Sint(M)$ and any monomorphism $g \colon V_J \to M$ in $\mod A$,
$g = fh$ for some $h \colon V_J \to X$.
\end{enumerate}
\end{prp}

\begin{proof}
(1) \implies (2$'$). This is trivial by Proposition \ref{prp:crt-r-int-appx}.

(2$'$) \implies (1). Assume the statement (2$'$).
It is enough to show the statement (2) in Proposition \ref{prp:crt-r-int-appx}.
Let $I \in \Sint(M)$ and $g \colon V_I \to M$ a monomorphism in $\mod A$.
Since $(\Sint(M), \le_S)$ is a finite poset, there exists some $J \in \max\Sint(M)$
such that $I \le_S J$.
Then by definition of the partial order $\le_S$, there exists a monomorphism
$g' \colon V_J \rightarrowtail M$ that is an extension of $g$ along $\si_{J,I}$.
Then by (2$'$), we have $g' = fh$ for some $h \colon V_J \to X$ as in the commutative diagram:
\begin{equation}
\label{eq:Sint-ext}
\begin{tikzcd}
	{V_J} & {V_I} \\
	X & M
	\arrow["h"', dashed, from=1-1, to=2-1]
	\arrow["{g'}", tail, from=1-1, to=2-2]
	\arrow["\si_{J,I}"', hook', from=1-2, to=1-1]
	\arrow["g", tail, from=1-2, to=2-2]
	\arrow["f"', from=2-1, to=2-2]
\end{tikzcd}.
\end{equation}
Thus (2) holds.
\end{proof}

\begin{rmk}
\label{rmk:improvements}
Corollaries \ref{cor:const-r-int-appx}, \ref{cor:form-min-r-appx}
are improved as well by replacing $\Sint(M)$ with $\max\Sint(M)$.
The dual remark works for Proposition \ref{prp:crt-r-int-appx}$'$ below.
These are used in Example \ref{exm:CL5}.
\end{rmk}

The notion dual to Definition \ref{dfn:order-S} is given as follows:

\begin{dfn}
Let $M \in \mod A$. Then we define a partial order $\le_F$ on $\Fint(M)$ as follows:
For any $I,\, J \in \Fint(M)$,
$I \le_F J$ if and only if $I$ is a down-set in $J$ such that 
any epimorphism $M \twoheadrightarrow V_I$ is lifted to an epimorphism $M \twoheadrightarrow V_J$
along $\pi_{I,J}$.
We set
\[
\max \Fint(M):= \{I \in \Fint(M) \mid I \text{ is maximal in } (\Fint(M), \le_F)\}.
\]
\end{dfn}

The statements dual to Propositions \ref{prp:crt-r-int-appx}, \ref{prp:smaller-Sint} and 
Corollaries \ref{cor:const-r-int-appx}, \ref{cor:form-min-r-appx} hold,
to which we refer as Proposition \ref{prp:crt-r-int-appx}$'$ (the combined version) and
Corollaries \ref{cor:const-r-int-appx}$'$, \ref{cor:form-min-r-appx}$'$,
respectively.
For example, Proposition \ref{prp:crt-r-int-appx}$'$ is stated as follows:

\newtheorem*{prp-a}{Proposition \ref{prp:crt-r-int-appx}$'$}
\begin{prp-a}
Let $f \colon X \to M$ be in $\mod A$.
Then the following are equivalent.
\begin{enumerate}
\item
$f$ is a left interval approximation of $M$.
\item
For any $I \in \Fint(M)$ and any epimorphism $g \colon M \to V_I$ in $\mod A$,
$g = hf$ for some $h \colon X \to V_I$.

\item[(2$'$)]
For any $I \in \max\Fint(M)$ and any epimorphism $g \colon M \to V_I$ in $\mod A$,
$g = hf$ for some $h \colon X \to V_I$.
\end{enumerate}
\end{prp-a}

The following is immediate from Propositions \ref{prp:crt-r-int-appx}, \ref{prp:smaller-Sint}
and \ref{prp:crt-r-int-appx}$'$, which can be used to check whether a given right $A$-module is
interval decomposable or not.

\begin{cor}
\label{cor:int-dec-2}
Let $M \in \mod A$. Then the following are equivalent.
\begin{enumerate}
\item
$M$ is interval decomposable;
\item
$M \iso \Ds_{I \in \Sint(M)}V_I^{\be^0_M(I)}$;
\item
$\dim_\k M = \sum_{I \in \Sint(M)} \be^0_M(I)\dim_\k V_I$;
\item[(2$'$)]
$M \iso \Ds_{I \in \Fint(M)}V_I^{\ovl{\be}^0_M(I)}$;
\item[(3$'$)]
$\dim_\k M = \sum_{I \in \Fint(M)} \ovl{\be}^0_M(I)\dim_\k V_I$;\text{ and}
\item
$M \iso \Ds_{I \in \Sint(M)\cap \Fint(M)}V_I^{\be^0_M(I)} \iso \Ds_{I \in \Sint(M)\cap \Fint(M)}V_I^{\ovl{\be}^0_M(I)}$.
\end{enumerate}

Moreover, the statements remain valid even if $\Sint$ and $\Fint$
are replaced with $\max\Sint$ and $\max\Fint$, respectively.
In particular,
if $V_I$ is a direct summand of $M$, then
$I \in \max\Sint(M) \cap \max\Fint(M)$.
\end{cor}

\begin{proof}
Consider a minimal right interval approximation $f_0$
and a minimal left interval approximation $g_0$ of $M$
to have the following exact sequences:
\begin{equation}
\label{eq:int-apprx-exact}
\Ds_{I \in \Sint(M)}V_I^{\be^0_M(I)} \ya{f_0} M \to 0,
\qquad
0 \to M \ya{g_o} \Ds_{I \in \Fint(M)}V_I^{\ovl{\be}^0_M(I)}.
\end{equation}
Denote by (f) the fact that  $f_0$ is an isomorphism,
and by (g) the fact that $g_0$ is an isomorphism.
Then note that (f) \equivalent (1) \equivalent (g).
It is obvious by the exact sequences above that
(f) \equivalent (2) \equivalent (3), and that
(g) \equivalent (2$'$) \equivalent (3$'$).
Thus (1), (2), (3), (2$'$), and (3$'$) are equivalent.

(1) \implies (4).
Assume (1). Then by the above, we have
\begin{equation}\label{eq:Sint-Fint}
\Ds_{I \in \Sint(M)}V_I^{\be^0_M(I)} \iso M \iso
\Ds_{I \in \Fint(M)}V_I^{\ovl{\be}^0_M(I)}.
\end{equation}
By the Krull--Schmidt theorem, 
there exists a unique subset $S$ of $\bbI$ and a map $a \colon S \to \bbZ_{> 0}$
such that $M \iso \Ds_{I \in S}V_I^{a(I)}$.
Hence \eqref{eq:Sint-Fint} shows that $S \subseteq \Sint(M) \cap \Fint(M)$ and that
for each $I \in S$, we have $\be^0_M(I) = a(I) = \ovl{\be}^0_M(I) \ne 0$,
which shows (4). 

(4) \implies (1). This is trivial.

By Remark \ref{rmk:improvements}, note that both $f_0$ and $g_0$ have the improved forms
that are obtained from \eqref{eq:int-apprx-exact} by replacing
$\Sint(M)$ and $\Fint(M)$ with $\max\Sint(M)$ and $\max\Fint(M)$,
respectively.
Then by the same argument, we see that the additional statements hold.
\end{proof}

\begin{rmk}
In the above, the statement (4) implies the following:
\begin{enumerate}
\item[(5)]
$\dim_\k M = \sum_{I \in \Sint(M)\cap \Fint(M)} \be^0_M(I)\dim_\k V_I$\\
\indent \hspace{32pt}$= \sum_{I \in \Sint(M)\cap \Fint(M)} \ovl{\be}^0_M(I)\dim_\k V_I$.
\end{enumerate}
However, it is not clear whether the converse holds or not because we do not know the existence of
an epimorphism $\Ds_{I \in \Sint(M)\cap \Fint(M)} V_I^{\be_M^0(I)} \to M$ or a monomorphism
$M \to \Ds_{I \in \Sint(M)\cap \Fint(M)} V_I^{\ovl{\be}_M^0(I)}$.
\end{rmk}

In the above corollary, the values of $\be^0_M(I)$ or $\ovl{\be}^0_M(I)$
can be computed by the formula \eqref{eq:beta0-formula} or its dual version.
For example, we have the following.

\begin{cor}
Let $M \in \mod A$. Then $M$ is interval decomposable if and only if
the following holds:
\begin{enumerate}
\item[(3$^{*})$]
$\dim_\k M = \sum_{I \in \max\Sint(M)} \dim_\k H_0(\calK_I(M)\down)\dim_\k V_I$.
\end{enumerate}
\end{cor}

\subsection{Some examples}
We define a bound quiver $(Q, \ro):= CL_n$\ ($n \ge 2$) to be the following quiver
with full commutativity relations (call it a {\em commutative ladder}), which is used
in the following examples:
\begin{equation}
\label{eq:CL}
\begin{tikzcd}
\ovl{1} & \ovl{2} & \cdots & \ovl{n}\\ 
1 & 2 & \cdots & n
\Ar{1-1}{1-2}{"\ovl{a}_1"}
\Ar{1-2}{1-3}{"\ovl{a}_2"}
\Ar{1-3}{1-4}{"\ovl{a}_{n-1}"}
%%%%%%%
\Ar{2-1}{2-2}{"a_1"'}
\Ar{2-2}{2-3}{"a_2"'}
\Ar{2-3}{2-4}{"a_{n-1}"'}
%%%%%%%
\Ar{2-1}{1-1}{"b_1"}
\Ar{2-2}{1-2}{"b_2"}
\Ar{2-3}{1-3}{"\cdots" description, phantom}
\Ar{2-4}{1-4}{"b_n"'}
\end{tikzcd}.
\end{equation}
In the first example (Example \ref{exm:CL3}), we have the whole AR-quiver that has only 29 vertices (only 2 of them are not interval modules),
and it is easy to control morphisms between indecomposables.
In the second example (Example \ref{exm:CL5}), since the algebra $A$ is not representation-finite,
we do not use the whole AR-quiver.
However, since the computation of almost split sequences starting from (or ending in)
interval modules are relatively easy as explained in \cite[Sect.\ 5]{ABENY},
we can make full use of them in the computation.
We will apply Corollary \ref{cor:form-min-r-appx}
(resp.\ \ref{cor:form-min-r-appx}$'$) to
compute minimal right (resp.\ left) interval approximations.

In examples, each indecomposable module $M$ is denoted by its dimension vector
 $\udim M$.
 For each interval (subquiver) $J$ of $Q$, we present $J$ by $\udim V_J$.

\begin{exm}
\label{exm:CL3}
Let $A:= \k(CL_3)$.
Then the AR-quiver of $A$ is as follows:
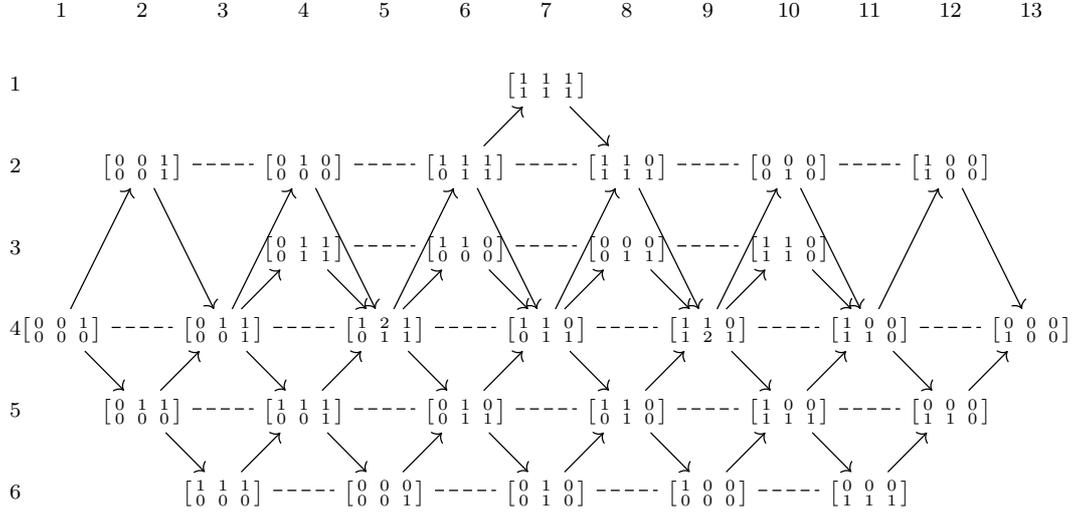
\[\tiny
\begin{tikzcd}[column sep=-7pt]
	& 1 & 2 & 3 & 4 & 5 & 6 & 7 & 8 & 9 & 10 & 11 & 12 & 13 \\
	1 &&&&&&& {\sbmat{1&1&1\\1&1&1}} \\
	2 && {\sbmat{0&0&1\\0&0&1}} && {\sbmat{0&1&0\\0&0&0}} && {\sbmat{1&1&1\\0&1&1}} && {\sbmat{1&1&0\\1&1&1}} && {\sbmat{0&0&0\\0&1&0}} && {\sbmat{1&0&0\\1&0&0}} \\
	3 &&&& {\sbmat{0&1&1\\0&1&1}} && {\sbmat{1&1&0\\0&0&0}} && {\sbmat{0&0&0\\0&1&1}} && {\sbmat{1&1&0\\1&1&0}} \\
	4 & {\sbmat{0&0&1\\0&0&0}} && {\sbmat{0&1&1\\0&0&1}} && {\sbmat{1&2&1\\0&1&1}} && {\sbmat{1&1&0\\0&1&1}} && {\sbmat{1&1&0\\1&2&1}} && {\sbmat{1&0&0\\1&1&0}} && {\sbmat{0&0&0\\1&0&0}} \\
	5 && {\sbmat{0&1&1\\0&0&0}} && {\sbmat{1&1&1\\0&0&1}} && {\sbmat{0&1&0\\0&1&1}} && {\sbmat{1&1&0\\0&1&0}} && {\sbmat{1&0&0\\1&1&1}} && {\sbmat{0&0&0\\1&1&0}} \\
	6 &&& {\sbmat{1&1&1\\0&0&0}} && {\sbmat{0&0&0\\0&0&1}} && {\sbmat{0&1&0\\0&1&0}} && {\sbmat{1&0&0\\0&0&0}} && {\sbmat{0&0&0\\1&1&1}}
	\arrow[from=5-2, to=6-3]
	\arrow[from=6-3, to=5-4]
	\arrow[from=5-4, to=3-5]
	\arrow[from=5-4, to=4-5]
	\arrow[from=5-4, to=6-5]
	\arrow[from=6-3, to=7-4]
	\arrow[from=7-4, to=6-5]
	\arrow[from=3-5, to=5-6]
	\arrow[from=4-5, to=5-6]
	\arrow[from=6-5, to=5-6]
	\arrow[dashed, no head, from=5-2, to=5-4]
	\arrow[from=5-6, to=3-7]
	\arrow[from=5-6, to=4-7]
	\arrow[from=5-6, to=6-7]
	\arrow[from=6-5, to=7-6]
	\arrow[from=7-6, to=6-7]
	\arrow[dashed, no head, from=6-3, to=6-5]
	\arrow[from=6-7, to=7-8]
	\arrow[from=4-7, to=5-8]
	\arrow[from=6-7, to=5-8]
	\arrow[from=3-7, to=5-8]
	\arrow[from=5-8, to=3-9]
	\arrow[from=5-8, to=4-9]
	\arrow[from=5-8, to=6-9]
	\arrow[from=7-8, to=6-9]
	\arrow[from=3-9, to=5-10]
	\arrow[from=4-9, to=5-10]
	\arrow[from=6-9, to=5-10]
	\arrow[from=6-9, to=7-10]
	\arrow[from=7-10, to=6-11]
	\arrow[from=5-10, to=6-11]
	\arrow[from=5-10, to=3-11]
	\arrow[from=5-10, to=4-11]
	\arrow[from=6-11, to=7-12]
	\arrow[from=3-11, to=5-12]
	\arrow[from=4-11, to=5-12]
	\arrow[from=6-11, to=5-12]
	\arrow[from=5-12, to=6-13]
	\arrow[from=7-12, to=6-13]
	\arrow[from=6-13, to=5-14]
	\arrow[dashed, no head, from=5-4, to=5-6]
	\arrow[dashed, no head, from=3-5, to=3-7]
	\arrow[dashed, no head, from=4-5, to=4-7]
	\arrow[dashed, no head, from=7-4, to=7-6]
	\arrow[dashed, no head, from=3-7, to=3-9]
	\arrow[dashed, no head, from=4-7, to=4-9]
	\arrow[dashed, no head, from=5-6, to=5-8]
	\arrow[dashed, no head, from=6-5, to=6-7]
	\arrow[dashed, no head, from=7-6, to=7-8]
	\arrow[dashed, no head, from=6-7, to=6-9]
	\arrow[dashed, no head, from=5-8, to=5-10]
	\arrow[dashed, no head, from=7-8, to=7-10]
	\arrow[dashed, no head, from=3-9, to=3-11]
	\arrow[dashed, no head, from=4-9, to=4-11]
	\arrow[dashed, no head, from=6-9, to=6-11]
	\arrow[dashed, no head, from=7-10, to=7-12]
	\arrow[dashed, no head, from=5-10, to=5-12]
	\arrow[dashed, no head, from=5-12, to=5-14]
	\arrow[dashed, no head, from=6-11, to=6-13]
	\arrow[from=2-8, to=3-9]
	\arrow[from=3-7, to=2-8]
	\arrow[dashed, no head, from=3-11, to=3-13]
	\arrow[from=3-13, to=5-14]
	\arrow[from=5-12, to=3-13]
	\arrow[from=5-2, to=3-3]
	\arrow[from=3-3, to=5-4]
	\arrow[dashed, no head, from=3-3, to=3-5]
\end{tikzcd}\]
We denote by $M_{i, j}$ the indecomposable module with the coordinate
$(i, j)$\ $(1\le i \le 6,\ 1 \le j \le 13)$. 
All these modules are interval modules except for $M_{4,5} = \sbmat{1&2&1\\0&1&1}$ and $M_{4,9} = \sbmat{1&1&0\\1&2&1}$.
Irreducible morphisms starting from $M_{i,j}$ are denoted by $a_{i,j}, a'_{i,j}, a''_{i,j}$ from the top, e.g., the morphism $M_{4,3} \to M_{5,4}$ is denoted by
$a''_{4,3}$.

(1)
We now compute the minimal interval Koszul coresolution $\calK\up(V_I)$ of
$V_I = M_{5,4}$ for the interval subquiver $I = \sbmat{1&1&1\\0&0&1}$.
First, the source map from $V_I$ gives rise to the almost split sequence
\[
0 \to M_{5,4} \ya{\mu_0} M_{4,5} \ds M_{6,5} \ya{\ep_0} M_{5,6} \to 0.
\]
Here, $E_I = M_{4,5} \ds M_{6,5}$.
The minimal left interval coresolution of $M_{4,5}$ is given by the
almost split sequence 
\begin{equation}
\label{eq:min-l-int-cores-M45}
0 \to M_{4,5} \ya{f'_0} M_{2,6} \ds M_{3,6} \ds M_{5,6} \ya{g_0} M_{4,7} \to 0
\end{equation}
starting from $M_{4,5}$ because both the middle term
and the end term are interval decomposable,
and that of $M_{6,5}$ is given by its identity.
Thus a minimal left interval approximation of $E_I$ is given by $f_0:= f'_0 \ds \id_{M_{6,5}}$,
which yields the short exact sequence
\[
0 \to M_{4,5} \ds  M_{6,5} \ya{f_0:= f'_0 \ds \id_{M_{6,5}}} (M_{2,6} \ds M_{3,6} \ds M_{5,6}) \ds  M_{6,5}
\ya{(g_0, 0)} M_{4,7} \to 0.
\]
Therefore, we have the commutative diagram
\[
\begin{tikzcd}[column sep=10pt]
0& M_{5,4} && X^1 && X^2 && \cdots\\
&M_{4,5}\ds M_{6,5}&& (M_{2,6}\ds M_{3,6} \ds M_{5,6})\ds M_{6,5} & C^1 && C^2 &&
\Ar{1-1}{1-2}{}
\Ar{1-2}{1-4}{"d^0"}
\Ar{1-2}{2-4}{"\et" '}
\Ar{1-4}{1-6}{"d^1"}
\Ar{1-6}{1-8}{"d^2"}
\Ar{1-2}{2-2}{"\mu_0"'}
\Ar{2-2}{2-4}{"f_0"'}
\Ar{2-4}{1-4}{"r_I"}
\Ar{1-4}{2-5}{"\mu_1"}
\Ar{2-5}{1-6}{"f_1"'}
\Ar{1-6}{2-7}{"\mu_2"'}
\Ar{2-7}{1-8}{"f_2"'}
\end{tikzcd}.
\]
Here, since $\mu_0 = \sbmat{a_{5,4}\\a'_{5,4}}$ and
$f'_0 = \sbmat{a_{4,5}\\a'_{4,5}\\a''_{4,5}}$,
 we have $\et = \sbmat{f_0 & 0\\0&\id}\sbmat{a_{5,4}\\a'_{5,4}}
 = \sbmat{f'_0 a_{5,4}\\a'_{5,4}}$, and hence
\[
\et = \sbmat{a_{4,5}a_{5,4}\\a'_{4,5}a_{5,4}\\a''_{4,5}a_{5,4}\\a'_{5,4}} = 
\sbmat{a_{4,5}a_{5,4}\\a'_{4,5}a_{5,4}\\-a_{6,5}a'_{5,4}\\a'_{5,4}} = 
\sbmat{1&0&0&0\\0&1&0&0\\0&0&1&-a_{6,5}\\0&0&0&1}\sbmat{a_{4,5}a_{5,4}\\a'_{4,5}a_{5,4}\\0\\a'_{5,4}},
\]
which is not left minimal.
By Lemma \ref{lem:to-r-min} (see also Remark \ref{rmk:r-min}(2)),
a left minimal version of $\et$ is given by
\[
d^0:= \sbmat{a_{4,5}a_{5,4}\\a'_{4,5}a_{5,4}\\a'_{5,4}} \colon
M_{5,4} \to \overbrace{M_{2,6} \ds M_{3,6} \ds M_{6,5}}^{X_1}.
\]
By the mesh relations\footnote{%%%%%%%%%%%%%%%%%%%%%%%%%%%%%
They are relations of the form $gf = 0$,
where $0 \to X \ya{f} Y \ya{g} Z \to 0$ is an almost split sequence.
}, we have
\begin{equation}
\label{eq:Im_in_Ker}
\sbmat{a'_{2,6}&a_{3,6}&-a_{5,6}a_{6,5}}
\sbmat{a_{4,5}a_{5,4}\\a'_{4,5}a_{5,4}\\a'_{5,4}}
= (a'_{2,6}a_{4,5}+a_{3,6}a'_{4,5}+a_{5,6}a''_{4,5})a_{5,4} = 0.
\end{equation}
On the other hand, it is easy to see that $d^0$ is a monomorphism and
the morphism
\[
\sbmat{a'_{2,6}&a_{3,6}&-a_{5,6}a_{6,5}} \colon
M_{2,6}\ds M_{3,6} \ds M_{6,5} \to M_{4,7}
\]
is an epimorphism by the almost split sequences in the AR-quiver.
Hence \eqref{eq:Im_in_Ker} shows that the morphism above is
a cokernel morphism of $d^0$ by counting the dimensions, and that
we may take $\mu_1:= \sbmat{a'_{2,6}&a_{3,6}&-a_{5,6}a_{6,5}}$
and $C^1:= M_{4,7}$.
Since $C^1$ is already interval decomposable, we see that the minimal Koszul coresolution $\calK\up(V_I)$
of $V_I$ is given by
\[
0 \to V_I \ya{d^0} \overbrace{M_{2,6}\oplus M_{3,6} \oplus M_{6,5}}^{X^1}
\ya{\mu_1} \overbrace{M_{4,7}}^{X^2} \to 0.
\]
In contrast, the pre-minimal interval Koszul coresolution $\bcalK\up(V_I)$ is given by
\[
0 \to V_I \ya{\et} M_{2,6}\oplus M_{3,6} \ds M_{5,6} \oplus M_{6,5}
\ya{\sbmat{a'_{2,6}&a_{3,6}&a_{5,6}&0\\0&0&\id&a_{6,5}}} M_{4,7}\ds M_{5,6} \to 0.
\]
For the module $M = M_{4,5}$, we compute the minimal interval Koszul complex $\calK_I(M)\down$.
From the AR-quiver, we immediately see that $\dm_A(X^1, M) = 0 = \dm_A(X^2, M)$,
which shows that
$\calK_I(M)\down = \dm_A(V_I, M) = \k a_{5,4}$
is a stalk complex concentrated in degree 0.
Hence we have $\be_M^i(I) = \de_{i,0}$ for all $i \ge 0$.
This coincides with the result obtained by looking at the interval resolution of $M$ that is given as the almost split sequence ending in $M$:
\[
0 \to M_{4,3} \to M_{2,4} \ds M_{3,4} \ds M_{5,4} \to M \to 0.
\]
Here, since $V_I = M_{5,4}$, we have $\be_M^i(I) = \de_{i,0}$ for all $i \ge 0$.
This also shows that $\de^\xi_M(I) = 1$
(see Definition \ref{dfn:int-app} and Corollary \ref{cor:int-appx}).

(2) We next compute the minimal interval Koszul coresolution $\calK\up(V_I)$ of
$V_I = M_{4,3}$ for the interval subquiver $I = \sbmat{0&1&1\\0&0&1}$.
The source map from $V_I$ gives rise to the almost split sequence
\begin{equation}
\label{eq:int-resol-M}
0 \to M_{4,3} \ya{\mu_0} M_{2,4} \ds M_{3,4} \ds M_{5,4} \ya{\ep_0} M_{4,5} \to 0.
\end{equation}
Here, $E_I = M_{2,4} \ds M_{3,4} \ds M_{5,4}$ is interval decomposable, and hence
the minimal left interval approximation of $E_I$ is given by the identity $\id_{E_I}$ of $E_I$,
and we see that $d^0 = \et = \mu_0$ and $C^1 = M_{4,5}$.
The minimal left interval coresolution of $M_{4,5}$ is given by the
almost split sequence \eqref{eq:min-l-int-cores-M45}.
As a consequence, the minimal interval Koszul coresolution $\calK\up(V_I)$ of $V_I$
is given as follows:
\[
0 \to M_{4,3} \ya{\mu_0} \overbrace{M_{2,4} \ds M_{3,4} \ds M_{5,4}}^{X^1} \ya{f_0 \ep_0}
 \overbrace{M_{2,6} \ds M_{3,6} \ds M_{5,6}}^{X^2}
\ya{g_0} \overbrace{M_{4,7}}^{X^3} \to 0.
\]
Then for the module $M = M_{4,5}$, the minimal interval Koszul complex $\calK_{I}(M)\down$ is given by
\[
\cdots\to 0 \to 0 \to \dm_A(X^1, M) \ya{\dm_A(\mu_0, M)} \dm_A(M_{4,3},M) \to 0
\]
because we see that $\dm_A(X^2, M) = 0 = \dm_A(X^3, M)$ by looking at the AR-quiver.
Since $\dm_A(\mu_0, M)$ is an epimorphism, we have $H_0(\calK_{I}(M)\down) = 0$,
and thus $\be_M^0(I) = 0$.
Now, since $\dm_A(X^1, M) = \k a_{2,4} \ds \k a_{3,4} \ds \k a_{5,4}$, each of its element
$f$ is uniquely written as $f = (r a_{2,4}, s a_{3,4}, t a_{5,4})$ 
for some $r, s, t \in \k$.
Then since
$\mu_0 = \sbmat{a_{4,3}\\a'_{4,3}\\ a''_{4,3}}$,
$\dm_A(\mu_0, M)(f) = 0$ if and only if
$r a_{2,4}a_{4,3} + s a_{3,4}a'_{4,3} + t a_{5,4}a''_{4,3} = 0$
if and only if $r = s = t$.
Hence
\[
H_1(\calK_{I}(M)\down) \iso \Ker \dm_A(\mu_0, M) = \k (a_{2,4}, a_{3,4}, a_{5,4}),
\]
and thus $\be_M^1(I) = 1$.
As a consequence, $\be_M^i(I) = \de_{i,1}$ for all $i \ge 0$,
which can also be read from the interval resolution \eqref{eq:int-resol-M} of $M$.
Moreover, we have $\de^\xi_M(I) = -1$.
\end{exm}

\begin{exm}
\label{exm:CL5}
Let $A:= \k(CL_5)$, $I:= \sbmat{0&0&0&1&1\\0&0&0&0&1}$, and $M \in \mod A$
be the indecomposable module
\begin{equation}\label{eq:str-M}
\begin{tikzcd}
0 & 0& \k & \k^2 & \k\\
0 & 0 & 0 & \k & \k
\Ar{1-1}{1-2}{}
\Ar{1-2}{1-3}{}
\Ar{1-3}{1-4}{"\sbmat{1\\1}"}
\Ar{1-4}{1-5}{"{(0,1)}"}
%%%%%
\Ar{2-1}{2-2}{}
\Ar{2-2}{2-3}{}
\Ar{2-3}{2-4}{}
\Ar{2-4}{2-5}{"1"}
%%%%%
\Ar{2-1}{1-1}{}
\Ar{2-2}{1-2}{}
\Ar{2-3}{1-3}{}
\Ar{2-4}{1-4}{"\sbmat{0\\1}"}
\Ar{2-5}{1-5}{"1"'}
\end{tikzcd}.
\end{equation}
Then $\udim M = \sbmat{0&0&1&2&1\\0&0&0&1&1}$.
We compute the interval Betti numbers $\be^i_M(I)$.
We first compute the minimal interval Koszul coresolution $\calK\up(V_I)$ of $V_I$.
The almost split sequence starting from $V_I$ is given as follows:
\begin{equation}
\label{eq:ass-I}
0 \to \sbmat{0&0&0&1&1\\0&0&0&0&1} \ya{\mu_0} \overbrace{\sbmat{0&0&0&1&1\\0&0&0&1&1} \ds\sbmat{0&0&0&1&0\\0&0&0&0&0} \ds \sbmat{0&0&1&1&1\\0&0&0&0&1}}^{X^1} \ya{\ep_0} M \to 0,
\end{equation}
where $\mu_0 = \dm^t(f_1,f_2,f_3)$ with $f_1, f_3$ canonical embeddings and $f_2$ a canonical
epimorphism.
Since the central term $E_I = X^1$ is already interval decomposable,
the minimal left interval approximation of $E_I$
is the identity $\id_{E_I}$ of $E_I$, and we have $d^0 = \et = \mu_0$.
Next, we can compute the minimal interval coresolution of $M$
by using Proposition \ref{prp:crt-r-int-appx}$'$, %and Proposition \ref{prp:smaller-Sint},
which is given as follows:
\begin{equation}
\label{eq:X2X3}
0 \to M \ya{\mu_1}
\overbrace{\sbmat{0&0&1&1&0\\0&0&0&0&0} \ds\sbmat{0&0&0&1&0\\0&0&0&1&1} \ds \sbmat{0&0&1&1&1\\0&0&0&1&1}}^{X^2}
\ya{d^2}
\overbrace{\sbmat{0&0&1&1&0\\0&0&0&1&1}}^{X^3} \to 0,
\end{equation}
where $\mu_1= \dm^t(g_1,g_2,g_3)$ with $g_1, g_2, g_3$ canonical epimorphisms having kernels
that are interval modules of dimension vectors
$\sbmat{0&0&0&1&1\\0&0&0&1&1}, \sbmat{0&0&1&1&1\\0&0&0&0&0}, \sbmat{0&0&0&1&0\\0&0&0&0&0}$,
respectively.
By combining these short exact sequences,
we obtain the minimal interval Koszul coresolution $\calK\up(V_I)$:
\[
0 \to V_I \ya{\mu_0} X^1 \ya{d^1} X^2 \ya{d^2} X^3 \to 0
\]
with $d^1:= \mu_1\ep_0$.
Hence the minimal interval Koszul complex $\calK_{I}(M)\down$ of $M$ at $I$ is given by
\[
0 \to \dm_A(X^3,M) \ya{\dm_A(d^2,M)} \dm_A(X^2,M) \ya{\dm_A(d^1,M)} \dm_A(X^1,M)
\ya{\dm_A(\mu_0,M)}\dm_A(V_I,M) \to 0.
\]
By looking at the structure \eqref{eq:str-M} of $M$ and those of indecomposable direct
summands of $X^2, X^3$ in \eqref{eq:X2X3},
we easily see that 
$\dm_A(X^3,M) = 0 = \dm_A(X^2,M)$.
Hence in particular, $\be^i_M(I) = 0$ for all $i \ge 2$, and $\Im \dm_A(d^1, M) = 0$.
We now compute $\Ker \dm_A(\mu_0,M)$.
Since $V_I$ is not a direct summand of $M$,
the almost split sequence \eqref{eq:ass-I} is mapped to the following exact sequence
by the functor $\dm_A(\blank, M)$:
\[
0 \to \End_A(M) \ya{\dm_A(\ep,M)} \dm_A(X^1, M) \ya{\dm_A(\mu_0,M)} \dm_A(V_I,M) \to 0.
\]
Thus $\Ker \dm_A(\mu_0,M) \iso \End_A(M) \iso \k$.
Hence
\[
\be^1_M(I) = \dim_\k H_1(\calK_{I}(M)\down) = \dim_\k \Ker \dm_A(\mu_0,M) = 1,
\]
and
\[
\be^0_M(I) = \dim_\k \Cok \dm_A(\mu_0, M) = 0.
\]
As a consequence, $\be^i_M(I) = \de_{i,1}$ for all $i \ge 0$, as the minimal interval
resolution \eqref{eq:ass-I} of $M$ suggests.
We also have $\de^\xi_M(I) = -1$.
\end{exm}

We now give an example, where $\La$ is isomorphic to the incidence algebra of
a finite lower semi-lattice, and compare the $\calI$-relative Koszul coresolution with
the formal Koszul coresolution in the introduction.

\begin{exm}
\label{exm:La=poset}
Let $\bfP$ be a poset whose Hasse quiver is given by the following quiver $Q$,
and set $A:= \k\bfP$:
$$
\begin{tikzcd}
&4\\
1 &2&3.
\Ar{2-1}{2-2}{}
\Ar{2-2}{2-3}{}
\Ar{2-2}{1-2}{}
\end{tikzcd}
$$
Consider a set $\calP:= \{V_I \mid I \in \bbI \setminus \{\{1\},\{3\},\{4\}\}\}$ 
of interval $A$-modules, where $\bbI$ is the set
of all intervals of the poset $\bfP$.
Regard $\calP$ as a full subcategory of $\mod A$ and set $\calI:= \add \calP$.
Then $\calP$ is isomorphic
to the incidence category of a finite lattice $\bfL$.
The AR-quiver of $A$ and the Hasse quiver of $\bfL$ is given as follows, respectively.
\[\begin{tikzcd}[column sep=14pt]
	&& {\sbmat{&1\\1&1&1}} && {\sbmat{&0\\0&1&0}} && {\sbmat{&0\\1&0&0}} \\
	{\sbmat{&1\\0&0&0}} && {\sbmat{&0\\0&1&1}} && {\sbmat{&1\\1&1&0}} \\
	& {\sbmat{&1\\0&1&1}} && {\sbmat{&1\\1&2&1}} && {\sbmat{&0\\1&1&0}} \\
	{\sbmat{&0\\0&0&1}} && {\sbmat{&1\\0&1&0}} && {\sbmat{&0\\1&1&1}}
	\arrow[dashed, no head, from=1-3, to=1-5]
	\arrow[from=1-3, to=3-4, "b_1" near start]
	\arrow[dashed, no head, from=1-5, to=1-7]
	\arrow[from=1-5, to=3-6, "d_1"]
	\arrow[dashed, no head, from=2-1, to=2-3]
	\arrow[from=2-1, to=3-2]
	\arrow[dashed, no head, from=2-3, to=2-5]
	\arrow[from=2-3, to=3-4, "b_2" ']
	\arrow[from=2-5, to=3-6, "d_2" ']
	\arrow[from=3-2, to=1-3, "a_1" near end]
	\arrow[from=3-2, to=2-3, "a_2" ']
	\arrow[dashed, no head, from=3-2, to=3-4]
	\arrow[from=3-2, to=4-3, "a_3" ']
	\arrow[from=3-4, to=1-5, "c_1" near end]
	\arrow[from=3-4, to=2-5, "c_2" ']
	\arrow[dashed, no head, from=3-4, to=3-6]
	\arrow[from=3-4, to=4-5, "c_3" ']
	\arrow[from=3-6, to=1-7]
	\arrow[from=4-1, to=3-2]
	\arrow[dashed, no head, from=4-1, to=4-3]
	\arrow[from=4-3, to=3-4, "b_3" ']
	\arrow[dashed, no head, from=4-3, to=4-5]
	\arrow[from=4-5, to=3-6, "d_3" ']
\end{tikzcd}\]
\[\begin{tikzcd}[column sep=22pt]
	&& {\sbmat{&1\\1&1&1}} && {\sbmat{&0\\0&1&0}} && {} \\
	{} && {\sbmat{&0\\0&1&1}} && {\sbmat{&1\\1&1&0}} \\
	& {\sbmat{&1\\0&1&1}} && {} && {\sbmat{&0\\1&1&0}} \\
	{} && {\sbmat{&1\\0&1&0}} && {\sbmat{&0\\1&1&1}}
	\arrow["{e_1}"{pos=0.2}, from=1-3, to=2-5]
	\arrow["{e_2}"{pos=0.8}, from=1-3, to=4-5, bend left=10pt]
	\arrow["{d_1}", from=1-5, to=3-6]
	\arrow["{f_1}"{pos=0.8}, from=2-3, to=1-5]
	\arrow["{f_2}"'{pos=0.8}, from=2-3, to=4-5]
	\arrow["{d_2}"', from=2-5, to=3-6]
	\arrow["{a_1}", from=3-2, to=1-3]
	\arrow["{a_2}"', from=3-2, to=2-3]
	\arrow["{a_3}"', from=3-2, to=4-3]
	\arrow["{g_1}"{pos=0.2}, from=4-3, to=1-5, bend left=10pt]
	\arrow["{g_2}"'{pos=0.3}, from=4-3, to=2-5]
	\arrow["{d_3}"', from=4-5, to=3-6]
\end{tikzcd},\]
where we set
$$
e_1:= c_2b_1,\, e_2:= -c_3b_1;\, f_1:= -c_1b_2,\, f_2:= c_3b_2;\, g_1:=c_1b_3,\, g_2:= -c_2b_3,
$$
which defines an isomorphism $\calP \to \k\bfL$.
For the interval $I = \sbmat{&1\\0&1&1}$, we compute
the minimal $\calI$-relative Koszul coresolution $\calK\up(V_I)$
(this coincides with the minimal interval Koszul coresolution in this case) and
the formal Koszul coresolution $C\up_{I}$ in Definition \ref{dfn:formalKos_coresol},
which is an $\calI$-relative Koszul coresolution by Corollary \ref{cor:weakKos_cores}.
They are given as follows, respectively:
$$
0 \to \sbmat{&1\\0&1&1} \ya{d^0}
\sbmat{&1\\1&1&1}\ds \sbmat{&0\\0&1&1} \ds \sbmat{&1\\0&1&0}
\ya{d^1}
\sbmat{&0\\0&1&0}\ds \sbmat{&1\\1&1&0} \ds \sbmat{&0\\1&1&1}
\ya{d^2}
\sbmat{&0\\1&1&0} \to 0;
$$
$$
0 \to \sbmat{&1\\0&1&1} \ya{d_C^0}
\sbmat{&1\\1&1&1}\ds \sbmat{&0\\0&1&1} \ds \sbmat{&1\\0&1&0}
\ya{d^1_C}
\sbmat{&0\\0&1&0}\ds \sbmat{&1\\1&1&0} \ds \sbmat{&0\\1&1&1}
\ya{d^2_C}
\sbmat{&0\\1&1&0} \to 0,
$$
where
$$
\begin{aligned}
d^0&:= \sbmat{a_1\\a_2\\a_3},\, d^1:= \sbmat{c_1\\c_2\\c_3}\sbmat{b_1\ b_2\ b_3},\, 
d^2:= \sbmat{d_1\ d_2\ d_3};\\
d_C^0&:= \sbmat{a_1\\a_2\\a_3},\, d^1_C:= \sbmat{0 &-f_1 & g_1\\-e_1 & 0 & g_2\\-e_2 & f_2 & 0},\, 
d^2_C:= \sbmat{d_1\ -d_2\ d_3}.
\end{aligned}
$$
Hence in this case, we have $\calK\up(V_I) \iso C\up_I$ as cocomplexes, not only in
the homotopy category.
\end{exm}

\section{Application to the interval approximation}
\label{sec:int-appx}

In this section, we give a direct application to the computation of
the \emph{compressed multiplicity} $c^\xi_M(I)$ and
the \emph{interval approximation} $\de^\xi_M(I)$ defined
in \cite[Definitions 5.1 and 5.6]{AENY-3}, respectively.

So far we used the language of algebras and modules over them
because this makes it easier to state and prove the statements.
However, in \cite{AENY-3}, we adopted the language of categories and modules over them
because we need a functor between categories, which is not interpreted
as a morphism between algebras preserving identity.
To apply our theorem to results in \cite{AENY-3},
we now give a bridge in the following.

\begin{rmk}
\label{rmk:alg-cat}
Let $(Q, \ro)$ be a bound quiver.
We here remark the relationship between the factor algebra $\k(Q, \ro)$ of the path algebra
$\k Q$ and the factor category $\k[Q,\ro]$ of the path category $\k[Q]$.
\begin{enumerate}
\item
A $\k$-linear category $\k[Q]$,
called the {\em path category} of $Q$ is defined as follows.
The set $\k[Q]_0$ of objects of $\k[Q]$ is equal to the vertex set $Q_0$ of $Q$,
and for any vertices $x, y \in Q_0$, $\k[Q](x,y) = e_x \k Q e_y$
under Convention \ref{cvn:left-right}
(also $\k[Q,\ro]$ uses the left-to-right notation).
The composite $p * q$ of $p \in \k[Q](x,y)$ and $q \in \k[Q](y,z)$ is defined to be
the product $pq$ in $\k Q$, and the identity
of each $x \in Q_0$ is the path $e_x$ of length 0 at $x$.
Then the matrix algebra $\Ds_{x,y \in Q_0}\k[Q](x,y)$ is isomorphic to $\k Q$.
\item
Each ideal $R$ of $\k Q$ corresponds the ideal $[R]:= (R(x,y):= e_x R e_y)_{x,y \in Q_0}$
of $\k[Q]$, and the bound quiver $(Q, \ro)$ defines a factor category
$\k[Q,\ro]:= \k[Q]/[\ro]$.
\item
Under Convention \ref{cvn:left-right},
the category $\mod \k[Q,\ro]$ of contravariant functors\footnote{%%%%%%%%%%%%
These are called right $\k[Q,\ro]$-modules.} from $\k[Q,\ro]$ to $\mod\k$ is 
isomorphic to the category $\rep(Q, \ro)$ of representations of $(Q, \ro)$,
and is equivalent to
the category $\mod \k(Q, \ro)$ of right $\k(Q, \ro)$-modules
.
\item
This point of view is needed when we consider a morphism $F$
from an algebra of the form $\k Q'$ with $Q'$ a finite quiver
to an algebra of the form $\k(Q, \ro)$ induced from
a quiver morphism $Q' \to U(\k[Q, \ro])$
that does not need to send the identity element of $\k Q'$ to
that of $\k(Q,\ro)$, where $U(\k[Q, \ro])$ is
the underlying quiver of $\k[Q, \ro]$,
namely the quiver obtained from
it by forgetting the map $\id \colon x \mapsto \id_x\ (x \in Q_0)$
and the composition of $\k[Q, \ro]$.

Even in that case, $F$ can be seen as a functor from the category $\k[Q']$ to
the category $\k[Q, \ro]$.
Therefore, in that situation, it is convenient to treat
the algebra $\k(Q, \ro)$ as the linear category $\k[Q, \ro]$.
This is done in this section.
\end{enumerate}
\end{rmk}

Throughout this section, we set $A:= \k(CL_n) = \k(Q, \ro)$  (see \eqref{eq:CL}), and
treat $A$ as the linear category $\k[Q, \ro]$ as explained in
Remark \ref{rmk:alg-cat}.
To state our application of Theorem \ref{thm:Kosz-Betti},
we recall necessary definitions and facts from \cite{AENY-3}.
The set $\bbI$ of all intervals $I$ in $Q$ are given as follows.
Since an interval $I$ is a full subquiver of $Q$, we denote it by its vertex set $I_0$.
For each $i, j \in \{1, \dots n\}$ with $i \le j$, we set
\[
\begin{aligned}{}
[i,j]&:= \{x \in \{1,\dots, n\} \mid i \le x \le j\} \subseteq Q_0,
\text{ and}\\
[\ovl{i},\ovl{j}]&:= \{\ovl{x} \in \{\ovl{1},\dots, \ovl{n}\} \mid i \le x \le j\} \subseteq Q_0.
\end{aligned}
\]
Then
by identifying $I$ with $I_0$ for all $I \in \bbI$,
\cite[Proposition 21]{ABENY} states that
\[
\bbI = \{ [\ovl{k}, \ovl{l}], [i,j], [\ovl{k}, \ovl{l}]\sqcup [i,j] \mid
1 \le k \le i \le l \le j \le n
\}.
\]
An interval $I$ with $I_0 = [\ovl{k},\ovl{l}] \sqcup [i,j]$ is
illustrated as follows:
\begin{equation}
\label{eq:fixed-interval}
\begin{tikzcd}
\ovl{k} & \cdots & \ovl{i} & \cdots & \ovl{l}\\
    &        & i & \cdots & l &\cdots & j
\Ar{1-1}{1-2}{"\ovl{a}_k"}
\Ar{1-2}{1-3}{"\ovl{a}_{i-1}"}
\Ar{1-3}{1-4}{"\ovl{a}_{i}"}
\Ar{1-4}{1-5}{"\ovl{a}_{l-1}"}
%%%%%%%%
\Ar{2-3}{2-4}{"a_i"'}
\Ar{2-4}{2-5}{"a_{l-1}"'}
\Ar{2-5}{2-6}{"a_{l}"'}
\Ar{2-6}{2-7}{"a_{j-1}"'}
%%%%%%%
\Ar{2-3}{1-3}{"b_i"}
\Ar{2-4}{1-4}{"\cdots", phantom}
\Ar{2-5}{1-5}{"b_l"'}
\end{tikzcd}
\end{equation}
Now let $Q'$ be the quiver:
\[
  \begin{tikzcd}
    & 2 && 4\\
    1 && 3 && 5
    \Ar{1-2}{2-1}{"\al_1"'}
    \Ar{1-2}{2-3}{"\al_2"}
    \Ar{1-4}{2-3}{"\al_3"'}
    \Ar{1-4}{2-5}{"\al_4"}
  \end{tikzcd}
\]
and set $B:= \k[Q']$ to be the path category of the quiver $Q'$.
We here regard $Q'_0$ as a totally ordered subset of the set of integers.
For any $i, j \in Q'_0$ with $i \le j$, we set $[i,j]:= \{x \in Q'_0 \mid i \le x \le j\}$.
Then the set $\bbI(Q')$ of interval subquivers of $Q'$ is given by
\[
\bbI(Q') = \{[i,j] \mid 1 \le i \le j \le 5\}
\]
that has 15 elements,
where as before we identify each interval subquiver with its vertex set.
By Gabriel's theorem the set $\calL:= \{V_I \mid I \in \bbI(Q')\}$ forms a complete set of
representatives of the isoclasses of indecomposable right $B$-modules.
Then each right $B$-module $M$ is uniquely decomposed as
\begin{equation}
\label{eq:d_M}
M = \Ds_{I \in \bbI(Q')}V_I^{d_M(V_I)}
\end{equation}
with $d_M(V_I) \ge 0$ by the Krull--Schmidt theorem, which defines a function
$d_M \colon \calL \to \bbZ_{\ge 0}$.

For each morphism $f$ in $\k[Q]$,
$[f]$ denotes the image of $f$ under the canonical functor $\k[Q] \to A$.
For each interval $I$, we define a quiver morphism $\xi_I \colon Q' \to U(A)$ as follows,
where $U(A)$ is the underlying quiver of the category $A$ (see Remark \ref{rmk:alg-cat}(4)).
\begin{enumerate}
\item In the case that $I_0 = [\ovl{k}, \ovl{l}] \sqcup [i,j]$
with $1 \le k \le i \le l \le j \le n$.
Define $\xi_I$ by the following table:
  \begin{center}
    \begin{tabular}{c||c|c|c|c|c|c|c|c|c}
      $x$ &  1 & 2 & 3 & 4 & 5 & $\al_1$ & $\al_2$ & $\al_3$ & $\al_4$\\
      \hline
      $\xi_I(x)$  \rule{0pt}{2.5ex}&$\ovl{l}$ & $\ovl{k}$ & $\ovl{i}$ & $i$ & $j$& $[\ovl{a}_k\cdots \ovl{a}_{l-1}]$ &$[\ovl{a}_k\cdots \ovl{a}_{i-1}]$ & $[b_i]$ & $[a_i\cdots a_{j-1}]$
    \end{tabular},
  \end{center}
  where if $k=i$, then $\ovl{a}_k\cdots \ovl{a}_{i-1}$ is replaced by $e_i$,
  (similar for the cases $i=l$ or $l=j$).

\item In the case that $I_0 = [i,j]$ with $0 \le i \le j \le n$.
  Define $\xi_I$ by the following table:
  \begin{center}
    \begin{tabular}{c||c|c|c|c|c|c|c|c|c}
      $x$ &  1 & 2 & 3 & 4 & 5 & $\al_1$ & $\al_2$ & $\al_3$ & $\al_4$\\
      \hline
      $\xi_I(x)$ & $j$ & $i$ & $i$ & $i$ & $j$& $[a_i\cdots a_{j-1}]$ &$[e_i]$ & $[e_i]$ & $[a_i\cdots a_{j-1}]$
    \end{tabular}.
  \end{center}

\item In the case that $I_0= [\ovl{k}, \ovl{l}]$ with $0 \le k \le l \le n$.
Define $\xi_I$ 
by the following table:
  \begin{center}
    \begin{tabular}{c||c|c|c|c|c|c|c|c|c}
      $x$ &  1 & 2 & 3 & 4 & 5 & $\al_1$ & $\al_2$ & $\al_3$ & $\al_4$\\
      \hline
      $\xi_I(x)$ \rule{0pt}{2.5ex}& $\ovl{l}$ & $\ovl{k}$ & $\ovl{k}$ & $\ovl{k}$ & $\ovl{l}$& $[\ovl{a}_k\cdots \ovl{a}_{l-1}]$ &$[e_k]$ & $[e_k]$ & $[\ovl{a}_k\cdots \ovl{a}_{l-1}]$
    \end{tabular}.
  \end{center}
\end{enumerate}

In the first case, $\xi_I$ is visualized as follows:
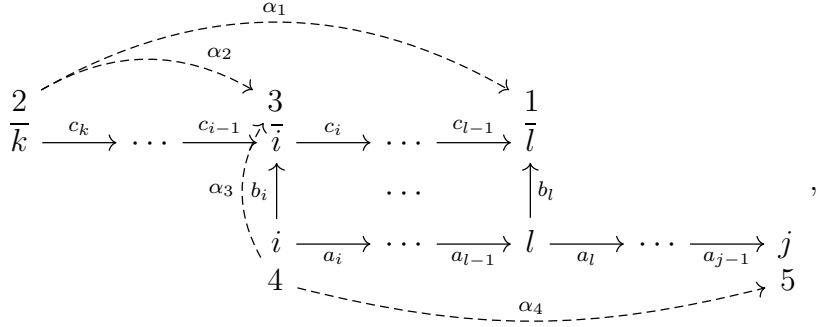
\[
\begin{tikzcd}
\ovl{k} & \cdots & \ovl{i} & \cdots & \ovl{l}\\
    &        & i & \cdots & l &\cdots & j
\Ar{1-1}{1-2}{"c_k"}
\Ar{1-2}{1-3}{"c_{i-1}"}
\Ar{1-3}{1-4}{"c_{i}"}
\Ar{1-4}{1-5}{"c_{l-1}"}
%%%%%%%%
\Ar{2-3}{2-4}{"a_i"'}
\Ar{2-4}{2-5}{"a_{l-1}"'}
\Ar{2-5}{2-6}{"a_{l}"'}
\Ar{2-6}{2-7}{"a_{j-1}"'}
%%%%%%%
\Ar{2-3}{1-3}{"b_i"}
\Ar{2-4}{1-4}{"\cdots", phantom}
\Ar{2-5}{1-5}{"b_l"'}
\end{tikzcd}
\hspace{-106.5mm}
\begin{tikzcd}[column sep=14mm, row sep=18mm]
2 &  & 3 &  & 1\\
    &        & 4 &  &  & & 5
\Ar{1-1}{1-5}{bend left, dashed,"\al_1"}
\Ar{1-1}{1-3}{bend left, dashed,"\al_2", near end}
\Ar{2-3}{1-3}{bend left, dashed,"\al_3"}
\Ar{2-3}{2-7}{bend right=5mm, dashed,"\al_4"}
\end{tikzcd},
\]
where each broken arrow $\al_i$ represent an arrow in the quiver $Q'$,
and the corresponding solid path represents its image $\xi_I(\al_i)$ in the quiver $U(A)$.

In any case, $\xi_I$ uniquely extends to a linear functor $F_I \colon B \to A$.
By using $F_I$, we regard $A$ to be the $B$-$A$-bimodule\footnote{%%%%%%%%
Note that by Convention \ref{cvn:left-right}, the bifunctor $(x,y) \mapsto A(x,y)$ is covariant in $x$
and contravariant in $y$.
Therefore, $(x,y) \mapsto A(F(x),y)$ is a $B$-$A$-bimodule.
} $\bi{A}BA = A(F_I(\cdot),\blank)$,
which gives us an adjoint pair
\[
  \begin{tikzcd}
    \mod A & \mod B
    \Ar{1-1}{1-2}{"R"', bend right=20}
    \Ar{1-2}{1-1}{"L"', bend right=20}
  \end{tikzcd},
\]
where
$L = L_{\xi_I}:=  \blank \otimes_B (\bi{A}BA)$ is a left adjoint to
$R = R_{\xi_I}:= \Hom_A(\bi{A}BA, \blank)$.
Recall that a right $A$-module $M$ is a contravariant functor $A \to \mod \k$
(Remark \ref{rmk:alg-cat}(3)), and then
$R(M) \iso M \circ F_I$.
For instance, $R(V_I) \iso V_{[1,5]}$.

\begin{dfn}[Compressed multiplicity]
  \label{def:compressedmult}
Let $\xi$ be the map sending $I \in \bbI$ to $\xi_I$ as defined above.
We define the \emph{compressed multiplicity with respect to $\xi$} of $V_I$ in $M$ as
\[
c_M^{\xi}(I):= d_{R_{\xi_I}(M)}(R_{\xi_I}(V_I)) = d_{R_{\xi_I}(M)}(V_{[1,5]})
\]
(see \eqref{eq:d_M}), which defines a map $c^\xi_M \in \bbR^\bbI$,
$I \mapsto c^\xi_M(I)\ (I \in \bbI)$.
\end{dfn}

We regard the set $\bbI$ as a poset by the partial order $\le$ defined by
$I \le J$ if and only if $I_0 \subseteq J_0$ for all $I, J \in \bbI$.
Consider the incidence algebra $R:= \bbR(\bbI\op) = \bbR(H(\bbI\op), \ro)$ (see Definition \ref{dfn:incidence} and Remark \ref{rmk:incidence}) of the opposite poset $\bbI\op$ of $\bbI$.
Then since $\bbI\op$ itself is an interval, we can consider the interval module
$V_{\bbI\op}$, the underlying vector space of which is equal to $\bbR^\bbI$.
By identifying these, $\bbR^\bbI$ can be seen as a right $R$-module.
It is well known that the zeta function $\ze:= \sum_{[I,J]\in \Seg(\bbI\op)}[I,J] \in R$
has an inverse $\mu$, called the M{\"o}bius function (see \cite{Rota} for details).

\begin{dfn}[Interval approximation]
\label{dfn:int-app}
We define the \emph{interval approximation}\footnote{%%%%%%%%%%%%%%%%%%%%%%%%%%%
This is also called {\em interval replacement}.}
$\de_M^{\xi}$ with respect to $\xi$ to be the
the M\"obius inversion of $c_M^{\xi}$, i.e., $\de^\xi_M:= c^\xi_M \mu$, the explicit form of
which is given by
\[
  \de_M^{\xi}(J):= \sum_{S\subseteq \cov(J)} (-1)^{\# S}c_M^{\xi}(\vee S)
\]
for all $J \in \bbI$.
\end{dfn}

We refer the reader to \cite{AENY-2, AENY-3, AGL} for the meanings and properties of 
the compressed multiplicity and the interval approximation.
The following is a main theorem in \cite[Sect.\ 3]{AENY-3}.
The second equality follows from the first by the general theory of M\"obius inversion.

\begin{thm}[Theorem 5.5 and Corollary 5.7 in \cite{AENY-3}]
\label{thm:AENY-3}
\label{cor:de-xi}
  Let $I \in \bbI$, and $M$ be in $\mod A$ with the minimal interval resolution \eqref{eq:min-r-app} with
$X_{r+1} = 0$ for some $r \ge 0$.
$($By \cite[Proposition 4.5]{AENY-3}, there always exists such an $r$.$)$
  Then we have
  \[
\begin{aligned}
c^{\xi}_M(I)&= \sum_{I \le J \in \bbI} \left(\sum_{i=0}^r (-1)^i \be_M^{i}(J)\right),
\text{ and}\\
    \de_M^{\xi}(I) &= \sum_{i=0}^r (-1)^i\be_M^{i}(I).
\end{aligned}
  \]
\end{thm}

We are now in a position to state our application of Theorem \ref{thm:Kosz-Betti}.
The following is a direct consequence of Theorems \ref{thm:Kosz-Betti}
and \ref{thm:AENY-3}.

\begin{cor}
\label{cor:int-appx}
Let $A = \k (CL_n)$ for some $n \ge 2$, $M \in \mod A$,
and $I \in \bbI$.
Suppose that $M$ has the minimal interval resolution \eqref{eq:min-r-app} with
$X_{r+1} = 0$ for some $r \ge 0$.
Then we have
\[
\begin{aligned}
c^\xi_M(I) &= \sum_{I \le J \in \bbI} \left(\sum_{i = 0}^r (-1)^i \dim_\k H_i(\calK_J(M)\down)\right), \text{ and}\\
\de^\xi_M(I) &= \sum_{i = 0}^r (-1)^i \dim_\k H_i(\calK_{I}(M)\down).
\end{aligned}
\]
\end{cor}

\end{document}